\documentclass{bmcart}
\startlocaldefs
\usepackage{amsmath} 
\usepackage{amsthm} 
\usepackage{amssymb}  
\usepackage{graphicx} 
\usepackage[utf8]{inputenc}
\usepackage{url}
\usepackage{epstopdf}
\usepackage{enumerate}

\newcommand{\cref}[1]{Chap.~\ref{#1}}  
\newcommand{\Ba}{B^{(a)}}
\newcommand{\Ra}{R^{(a)}}

\newcommand{\Va}{\mathbf{V}^{(a)}}
\newcommand{\Pci}{P_{c,i}}
\newcommand{\Scid}{S_{c,i}^{(d)}}
\newcommand{\Scia}{S_{c,i}^{(a)}}
\newcommand{\Xcia}{X_{c,i}^{(a)}}
\endlocaldefs
\begin{document}
\begin{frontmatter}
\begin{fmbox}
\dochead{Research}
\title{Prediction and Optimal Scheduling of Advertisements in Linear Television}

\author[addressref={aff1},
   email={panaggio@rose-hulman.edu},  
    corref={aff1},  
]{\inits{MJ}\fnm{Mark J} \snm{Panaggio}}
\author[addressref={aff8},
   email={pakwing@udel.edu },   
]{\inits{P-W}\fnm{Pak-Wing} \snm{Fok}}
\author[addressref={aff2},
   email={gbhatt@tnstate.edu},  
]{\inits{GS}\fnm{Ghan S} \snm{Bhatt}}
\author[addressref={aff3},
   email={burhoe@math.umass.edu }, 
]{\inits{S}\fnm{Simon} \snm{Burhoe}}
\author[addressref={aff4},
   email={cappsm@rams.colostate.edu},   
]{\inits{M}\fnm{Michael} \snm{Capps}}
\author[addressref={aff5},
   email={cedholm2@math.unl.edu },   
]{\inits{CJ}\fnm{Christina J} \snm{Edholm}}
\author[addressref={aff6},
   email={fadoua@temple.edu},   
]{\inits{F}\fnm{Fadoua} \snm{El Moustaid}}
\author[addressref={aff4},
   email={emerson@math.colostate.edu },   
]{\inits{T}\fnm{Tegan} \snm{Emerson}}
\author[addressref={aff6},
   email={starlenaq@gmail.com },   
]{\inits{S-L}\fnm{Star-Lena} \snm{Estock}}
\author[addressref={aff9},
   email={ngold5@my.yorku.ca},   
]{\inits{N}\fnm{Nathan} \snm{Gold}}
\author[addressref={aff10},
   email={rghalabi@math.ucdavis.edu},   
]{\inits{R}\fnm{Ryan} \snm{Halabi}}
\author[addressref={aff8},
   email={mhouser@udel.edu},   
]{\inits{M}\fnm{Madelyn} \snm{Houser}}
\author[addressref={aff12},
   email={kramep@rpi.edu },   
]{\inits{PR}\fnm{Peter R} \snm{Kramer}}
\author[addressref={aff13},
   email={waynelee1217@gmail.com },   
]{\inits{H-W}\fnm{Hsuan-Wei} \snm{Lee}}
\author[addressref={aff14},
   email={qli@fisk.edu },   
]{\inits{Q}\fnm{Qingxia} \snm{Li}}
\author[addressref={aff8},
   email={weiqiang@udel.edu },   
]{\inits{W}\fnm{Weiqiang} \snm{Li}}
\author[addressref={aff12},
   email={lud@rpi.edu  },   
]{\inits{D}\fnm{Dan} \snm{Lu}}
\author[addressref={aff12},
   email={qianyuzhou1@gmail.com },   
]{\inits{Y}\fnm{Yuzhou} \snm{Qian}}
\author[addressref={aff8},
   email={rossi@math.udel.edu },   
]{\inits{LF}\fnm{Louis F} \snm{Rossi}}
\author[addressref={aff20},
   email={dshutt@mines.edu },   
]{\inits{D}\fnm{Deborah} \snm{Shutt}}
\author[addressref={aff21},
   email={chuqiaoyang2013@u.northwestern.edu },   
]{\inits{C}\fnm{Vicky Chuqiao} \snm{Yang}}
\author[addressref={aff8},
   email={yxzhou@udel.edu},   
]{\inits{Y}\fnm{Yingxiang} \snm{Zhou}}
\address[id=aff1]{  \orgname{Mathematics Department, Rose-Hulman Institute of Technology }, 
  \city{Terre Haute, IN 47803},                              
  \cny{USA}                                    
}
\address[id=aff8]{  \orgname{Department of Mathematical Sciences, University of Delaware}, 
  \city{Newark, DE 19716},                              
  \cny{USA}                                    
}
\address[id=aff2]{  \orgname{Department of Mathematical Sciences, Tennessee State University}, 
  \city{Nashville, TN 37209},                              
  \cny{USA}                                    
}
\address[id=aff3]{  \orgname{Department of Mathematics and Statistics, University of Massachusetts Amherst}, 
  \city{Amherst, MA 01003},                              
  \cny{USA}                                    
}
\address[id=aff4]{  \orgname{Department of Mathematics, Colorado State University }, 
  \city{Fort Collins, CO 80523},                              
  \cny{USA}                                    
}
\address[id=aff5]{  \orgname{Department of Mathematics, University of Nebraska-Lincoln}, 
  \city{Lincoln, NE 68588},                              
  \cny{USA}                                    
}
\address[id=aff6]{  \orgname{Department of Mathematics, Temple University}, 
  \city{Philadelphia, PA 19122},                              
  \cny{USA}                                    
}
\address[id=aff9]{  \orgname{Department of Mathematics and Statistics, York University}, 
  \city{Toronto, ON M3J 1P3},                              
  \cny{Canada}                                    
}
\address[id=aff10]{  \orgname{Department of Mathematics, University of California, Davis}, 
  \city{Davis, CA 95616},                              
  \cny{USA}                                    
}
\address[id=aff12]{  \orgname{Department of Mathematical Sciences, Rensselaer Polytechnic Institute}, 
  \city{Troy, NY 12180},                              
  \cny{USA}                                    
}
\address[id=aff13]{  \orgname{Department of Mathematics, University of North Carolina at Chapel Hill }, 
  \city{Chapel Hill, NC 27599},                              
  \cny{USA}                                    
}
\address[id=aff14]{  \orgname{Department of Mathematics and Computer Science, Fisk University}, 
  \city{Nashville, TN 37208},                              
  \cny{USA}                                    
}
\address[id=aff20]{  \orgname{Applied Mathematics and Statistics Department, Colorado School of Mines}, 
  \city{Golden, CO 80401},                              
  \cny{USA}                                    
}
\address[id=aff21]{  \orgname{Department of Engineering Sciences and Applied Mathematics, Northwestern University}, 
  \city{Evanston, IL 60208},                              
  \cny{USA}                                    
}
\end{fmbox}

\begin{abstractbox}

\begin{abstract}
Advertising is a crucial component of marketing and an important way for companies to raise awareness of goods and services in the marketplace. Advertising campaigns are designed to convey a marketing image or message to an audience of potential consumers and television commercials can be an effective way of transmitting these messages to a large audience. In order to meet the requirements for a typical advertising order, television content providers must provide advertisers with a predetermined number of ``impressions'' in the target demographic. However, because the number of impressions for a given program is not known a priori and because there are a limited number of time slots available for commercials, scheduling advertisements efficiently can be a challenging computational problem.  In this case study, we compare a variety of methods for estimating future viewership patterns in a target demographic from past data. We also present a method for using those predictions to generate an optimal advertising schedule that satisfies campaign requirements while maximizing advertising revenue.
\end{abstract}




\end{abstractbox}
\end{frontmatter}

\section{Background}
The use of advertisements as a means for marketing consumer goods has been common practice for centuries. Mass distribution of advertisements through newspapers was first made possible by the printing press, but it was the advent of the radio and the television in the twentieth century that ultimately revolutionized advertising by allowing companies to transmit marketing messages into millions of homes around the world simultaneously \cite{obarr2010}. Today, advertising is an over \$500 trillion dollar global industry, and although advertising through digital media is growing rapidly, television remains the primary advertising medium with total television advertising expenditures making up approximately 40\% of the worldwide total \cite{McKinsey2014}.  

Advertising on ``linear'' (traditional live, not on-demand) television typically consists of an arrangement between content providers/programmers (TV networks such as ABC, NBC and Fox or cable operators such as Comcast or Cox) and advertising agencies in which the networks/operators are paid to run commercials in order to reach a desired audience. This audience is typically specified in demographic or psychographic terms, such as ``women 18 - 54'', or ``people concerned with health and fitness.'' A campaign's marketing target can be quantified by three measures: \textit{impressions}, which is the total number of times the message or ad is seen by a member of the target audience, \textit{reach}, which is the number of unique members of the target groups exposed to the ad, and \textit{frequency}, which is the average number of times the ad is viewed by each member of the target group that is reached by the ad. When a content provider agrees to fill an advertising order, they commit to running the commercial as many times as necessary until the desired number of impressions (and possibly reach and frequency) have been obtained.  Since the available commercial time in a given time frame is limited and each order that a content provider is able to fill provides additional revenue, it is of interest to meet each impression target in such a way that leaves broadcast time available for additional orders. Therefore, before an order can be accepted, content programmers must assess whether the number of impressions can be achieved in an acceptable time-frame and whether the budget for the order is large enough to warrant using broadcast time to provide those impressions. 

In order to make this determination, it is necessary to be able to generate a schedule that satisfies the order constraints in an efficient way. However, the generation of optimal advertising schedules poses a number of challenges.  First of all, one must know the viewership demographics of each television program.  Unfortunately, this is not known in advance and can only be estimated from past viewership data.  Secondly, depending on the content provider and time-frame, there can be a large number of possible orders and available commercial slots along with a large number constraints on what schedules are acceptable.  This makes finding an optimal schedule a computationally intensive problem that cannot usually be solved by hand. For this reason, naive approaches to scheduling lead to wasted resources and disenchanted audiences when ads fail to reach the interested consumers efficiently and must be aired repeatedly in order to meet impression targets. 

Here, we address two aspects of this interesting mathematical problem. First, using data modeled statistically mimic real TV viewership behavior as reported, for example by  The Nielsen Company \cite{Nielsen}, we explore a number of methods for predicting the number of impressions of future programming (Section \ref{sec:predict}) . These methods make use of spectral analysis, machine learning, Kalman filtering and distance scores.  Second, we demonstrate that when combined with advertising orders, these predictions can be used to formulate a nonlinear optimization problem that can be solved using standard integer programming techniques (Section \ref{sec:BIP}).  Finally, we outline a method for extending earlier results to account for the reach and frequency of an advertisement (Section \ref{sec:reach}).  This work was sponsored by clypd Inc.~and made possible by the 2015 Mathematical Problems in Industry Workshop.

\section{Predicting Number of Impressions}\label{sec:predict}
\begin{figure}[htbp]
\begin{center}
\includegraphics[width = 0.7\textwidth,clip=true,trim=0cm 0cm 0cm 0cm]{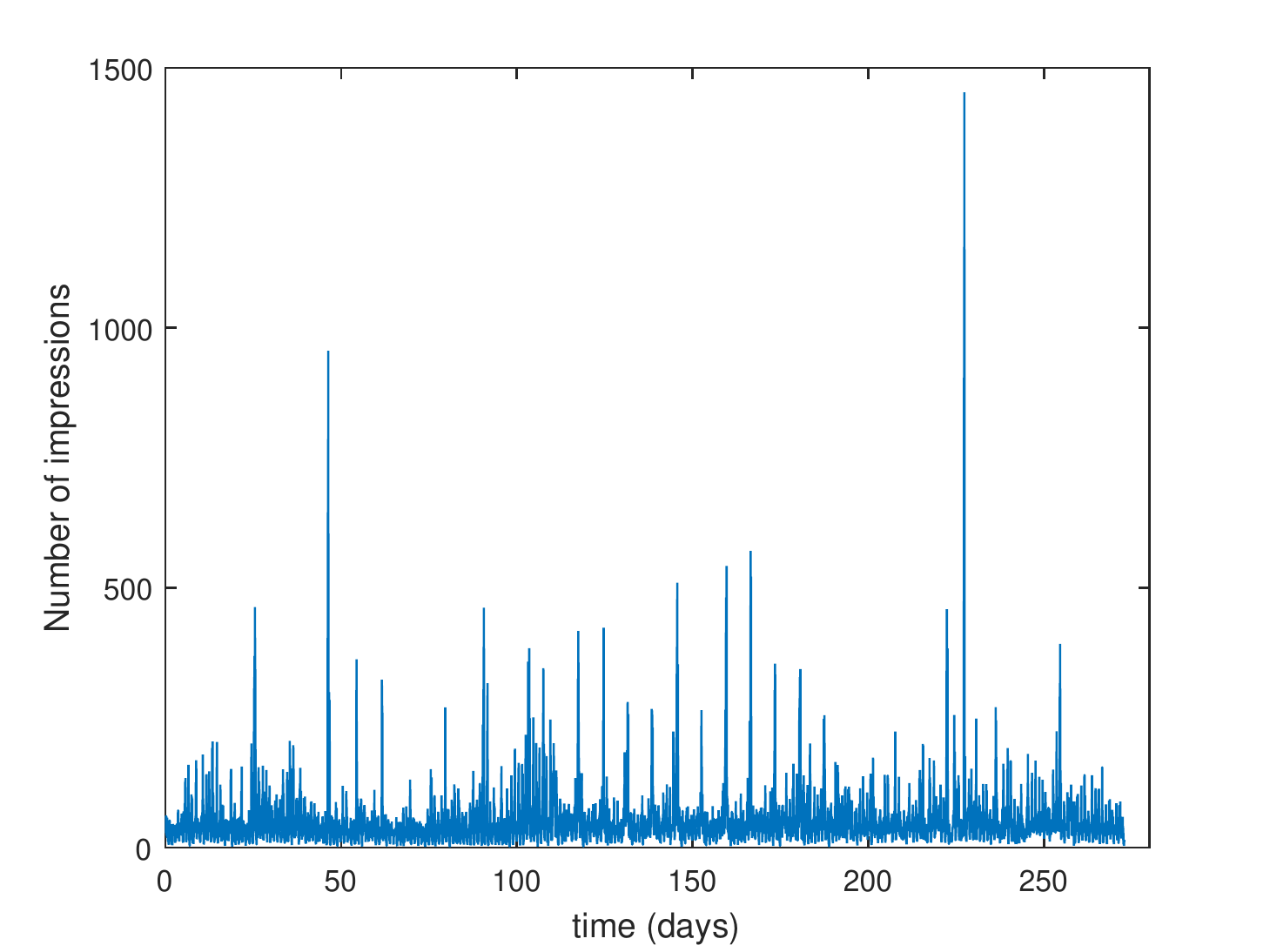}
\end{center}
\caption{Impressions per hour. Number of impressions for each hour over the course of 273 days starting
from October 23, 2014 on one particular channel.}
\label{fig:0}
\end{figure}

In order to generate predictions for the viewership and demographics of future programs, we used simulated past program viewership data provided by clypd Inc. Since precise data on the viewership of commercials was not available, we assume that program viewership data is representative of the viewership of advertisements as well.  Figure \ref{fig:0} shows a time series for the number of impressions on one particular channel over a period of about
273 days. Qualitatively, the signal is noisy with large spikes in the viewership. 

\subsection{Spectral Analysis}
\label{sec:spectral}
Although the data appears noisy, there are clear
periodic trends in the data. These are mostly likely driven by the periodic nature
of the channel programming. In order to identify these trends, we assume that the number of impressions, $S(t)$ can be decomposed
into a deterministic, periodic part $P(t)$ and a stochastic part $\eta(t)$,
\begin{equation}
S(t) = P(t) + \eta(t).
\label{eq:0} 
\end{equation}
We attempt to filter the signal to remove $\eta(t)$, leaving behind $P(t)$ by first filling in missing data using linear interpolation and then performing a
Fourier transform on the data and taking only the dominant modes in the power spectrum.

\begin{figure}[htbp]
\centering
\includegraphics[width=\textwidth,clip=true,trim=3cm 1cm 2cm 1cm]{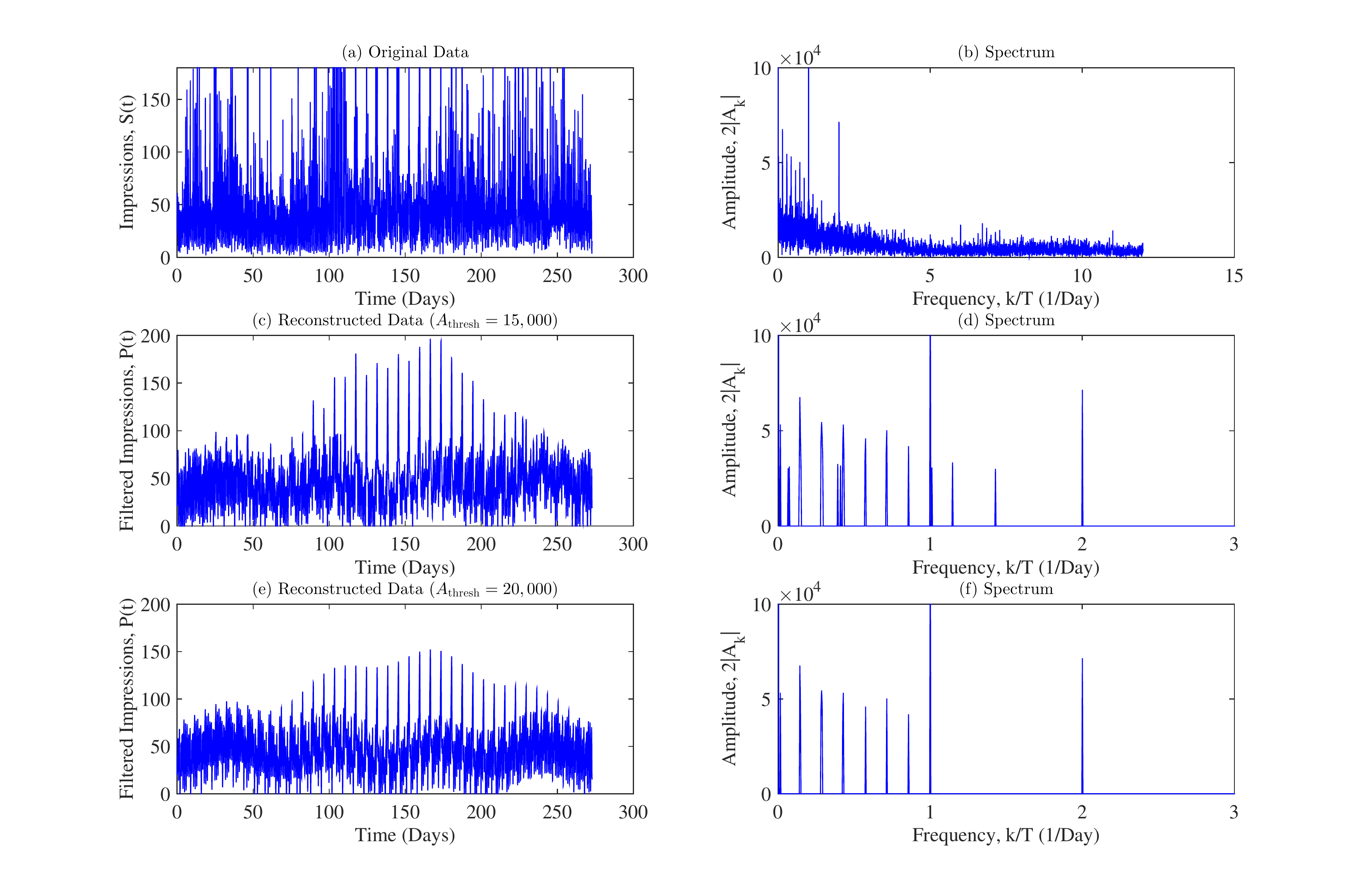}
\caption{Total viewership. (a,b): unfiltered data and power spectrum (frequency has units of day$^{-1}$).
(c,d): filtered signal and power spectrum with $A_{\textrm{thresh}} = 15,000$.
(e,f): filtered signal and power spectrum with $A_{\textrm{thresh}} = 20,000$.
}
\label{fig:1}
\end{figure}
\begin{figure}[htbp]
\centering
\includegraphics[width=\textwidth,clip=true,trim=3cm 1cm 2cm 1cm]{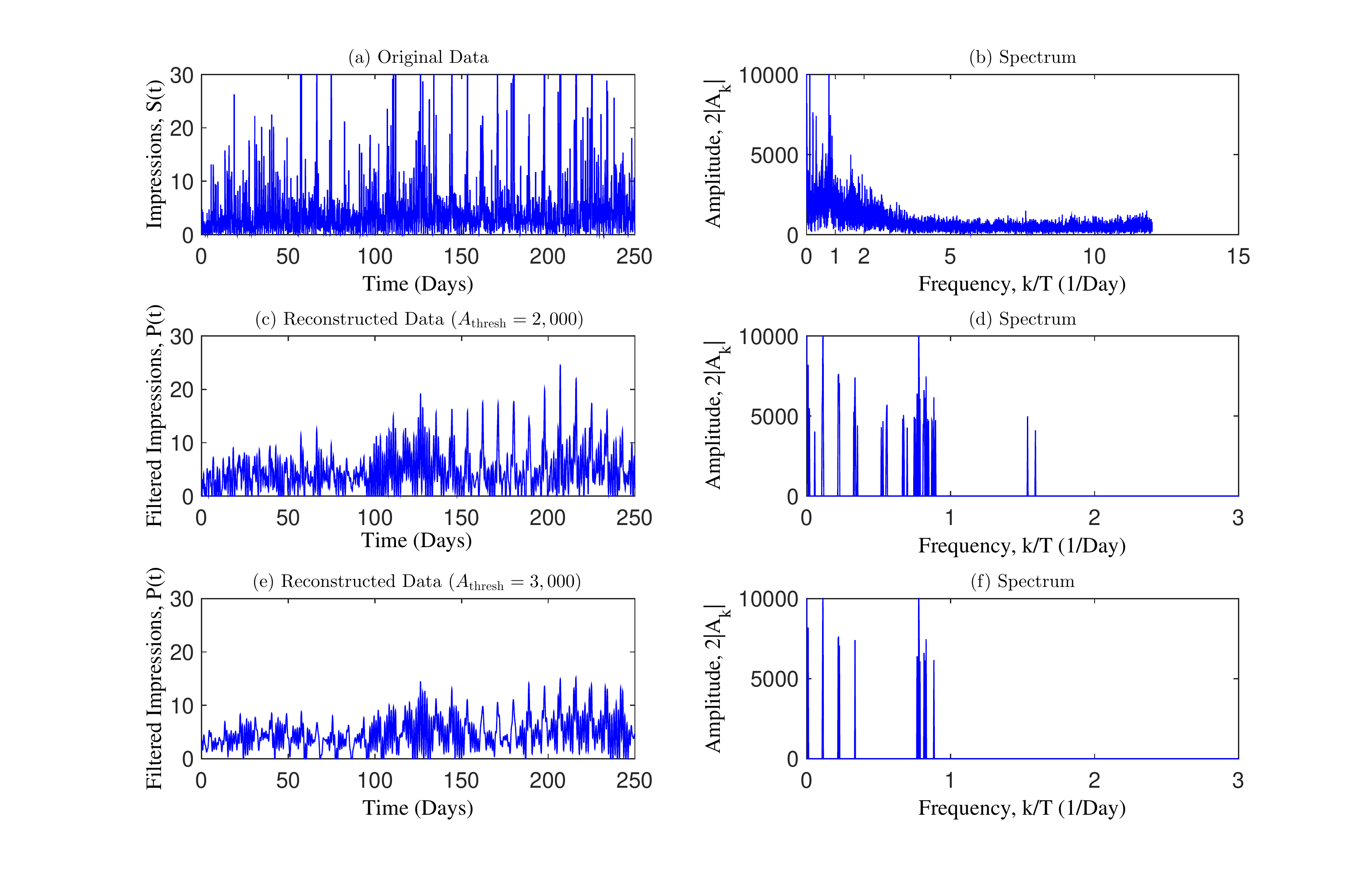}
\caption{Viewership of males 65 and older (a,b): unfiltered data and power spectrum (frequency has units of day$^{-1}$).
(c,d): filtered signal and power spectrum with $A_{\textrm{thresh}} = 2,000$.
(e,f): filtered signal and power spectrum with $A_{\textrm{thresh}} = 3,000$.
}
\label{fig:2}
\end{figure}

For this filtering scheme, we used Matlab's \texttt{fft} and \texttt{ifft} algorithms
to compute the power spectrum and then
removed all frequencies with amplitude less than
a given threshold $A_{\textrm{thresh}}$.
We write the full signal as
\begin{equation}
S(t) = A_0 + \sum_{j=1}^N (A_j e^{i k_j t/T} + A_j^* e^{-i k_j t/T} ),
\label{eqn:SS}
\end{equation}
where $t$ is the time (in hours), $N$ is the number of frequencies in the transform, $T$ is the total duration of the signal,  $A_j$ and $A_j^*$ are complex amplitudes (conjugates of each other), and $k_j$ are dimensionless frequencies. For our data set, we have $N = 3275$, and the total duration of the signal is $T = 6551$ hours $\approx 273$ days. This corresponds to 6551 data points $S(0),\ldots,S(6550)$ and 6551 Fourier coefficients (distinguishing between $A_j$ and $A_j^*$) and therefore, representation (\ref{eqn:SS}) exactly reproduces $S(t)$.

We decompose the signal $S(t)$ as follows 
\begin{equation}
S(t) = \underbrace{A_0 + \sum_{j: |A_j| > A_{\textrm{thresh}}} (A_j e^{i k_j t/T} + A_j^* e^{-i k_j t/T} )}_{\textrm{filtered signal}, ~P(t)}
+ \underbrace{\sum_{j: |A_j| \leq A_{\textrm{thresh}}} (A_j e^{i k_j t/T} + A_j^* e^{-i k_j t/T} )}_{
\textrm{noise},~\eta(t)}
\end{equation}
where the {\it signal} is defined as
\begin{equation}
P(t) = A_0 + \sum_{j: |A_j| > A_{\textrm{thresh}}} (A_j e^{i k_j t/T} + A_j^* e^{-i k_j t/T} ).
\end{equation}
while the {\it noise} is defined as
\begin{equation}
\eta(t) = \sum_{j: |A_j| \leq A_{\textrm{thresh}}} (A_j e^{i k_j t/T} + A_j^* e^{-i k_j t/T} ).
\end{equation}
Therefore, noise can be eliminated by removing all frequency  modes $j$ such that $|A_j| \leq A_{\textrm{thresh}}$. In other words, which Fourier modes are considered part of the signal and which are part of the noise depends solely
on the cut-off $A_{\rm thresh}$.

\subsubsection{Analysis of the signal}
The resulting signal is displayed in Figures \ref{fig:1} (for the general viewership) and \ref{fig:2}
(for males 65 and older). Power spectra are shown for rescaled frequencies $k_j/T$ which
have units of inverse time.
In Figure \ref{fig:1}, we see that
there are ``spikes'' at $k/T = i/7$ per day for $i=1,\ldots,7$. There are also
spikes at $k/T=0$ and $k/T = 2$ per day. Broadly speaking, the power spectra reveal that the predictable portion of the 
signal contains nine dominant frequencies/periods for viewer behavior
\begin{enumerate}
\item[1.] The zero mode in which the television is always on or off, regardless of the time ($k/T=0$ per day)
\item[2.] Viewing patterns that repeat weekly ($k/T=1/7$ per day)
\item[3.] Viewing patterns that repeat twice per week ($k/T=2/7$ per day)
\item[4.] Viewing patterns that repeat three times per week ($k/T=3/7$ per day)
\item[5.] Viewing patterns that repeat four times per week ($k/T=4/7$ per day)
\item[6.] Viewing patterns that repeat five times per week ($k/T=5/7$ per day)
\item[7.] Viewing patterns that repeat six times per week ($k/T=6/7$ per day)
\item[8.] Viewing patterns that repeat daily ($k/T=1$ per day)
\item[9.] Viewing patterns that repeat twice daily ($k/T=2$ per day)
\end{enumerate}
These behaviors do not necessarily refer to the same members of the viewership. 
The frequencies above also appeared in different demographic groups. However, males and females
65 and older did not exhibit these frequencies (see Fig. \ref{fig:2}) suggesting that
older members of the population have qualitatively different viewing habits.
We should note, however, that
upon taking different 7-week subsets of the data, some of the spikes at $k/T=1/7, 2/7,3/7,\ldots,6/7$ no longer appear. That is, the $k/T=0$, $k/T=1$ and $k/T=2$ per day modes are the most robust.

\begin{figure}[!htb]
\centering 
\includegraphics[width = 0.9\textwidth]{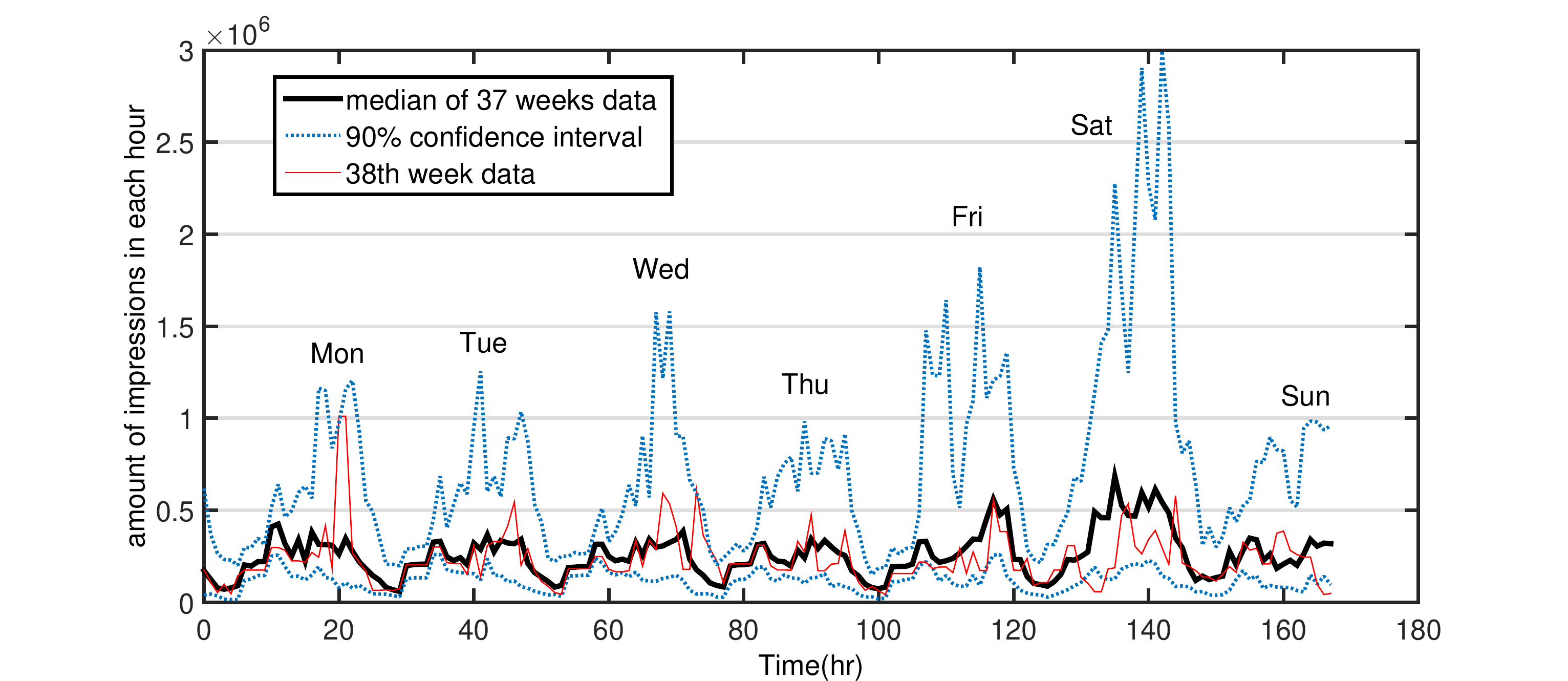}
\caption{Predicted week-long viewership trend. The signal (for channel 7287) exhibits
clear periodicity with dominant frequencies $k/T = \frac{1}{7},\frac{2}{7},\ldots,\frac{7}{7}$ per day
and $k/T=2$ per day.}
\label{v2}
\end{figure}

In summary, viewership patterns are periodic with both daily and weekly frequency components.
The weekly pattern is shown as a solid black curve in Figure \ref{v2}. This figure indicates the median number of impressions over a typical week along with confidence intervals. 
The distribution of impressions at a given time in the week is found from the first 38 weeks of
data. 90\% confidence intervals (dashed light blue) are calculated by taking the 5th and 95th percentiles of
the distribution. We see that during the evening of each day, there is a rise in the viewership.
Saturday seems to be the hardest to predict since it has the highest variability in impressions.
The red curve shows the viewership over the 39th week, which mostly falls within the confidence
intervals.

\subsubsection{Analysis of the noise}
\begin{figure}[htbp]
\centering
\includegraphics[width=0.7\textwidth]{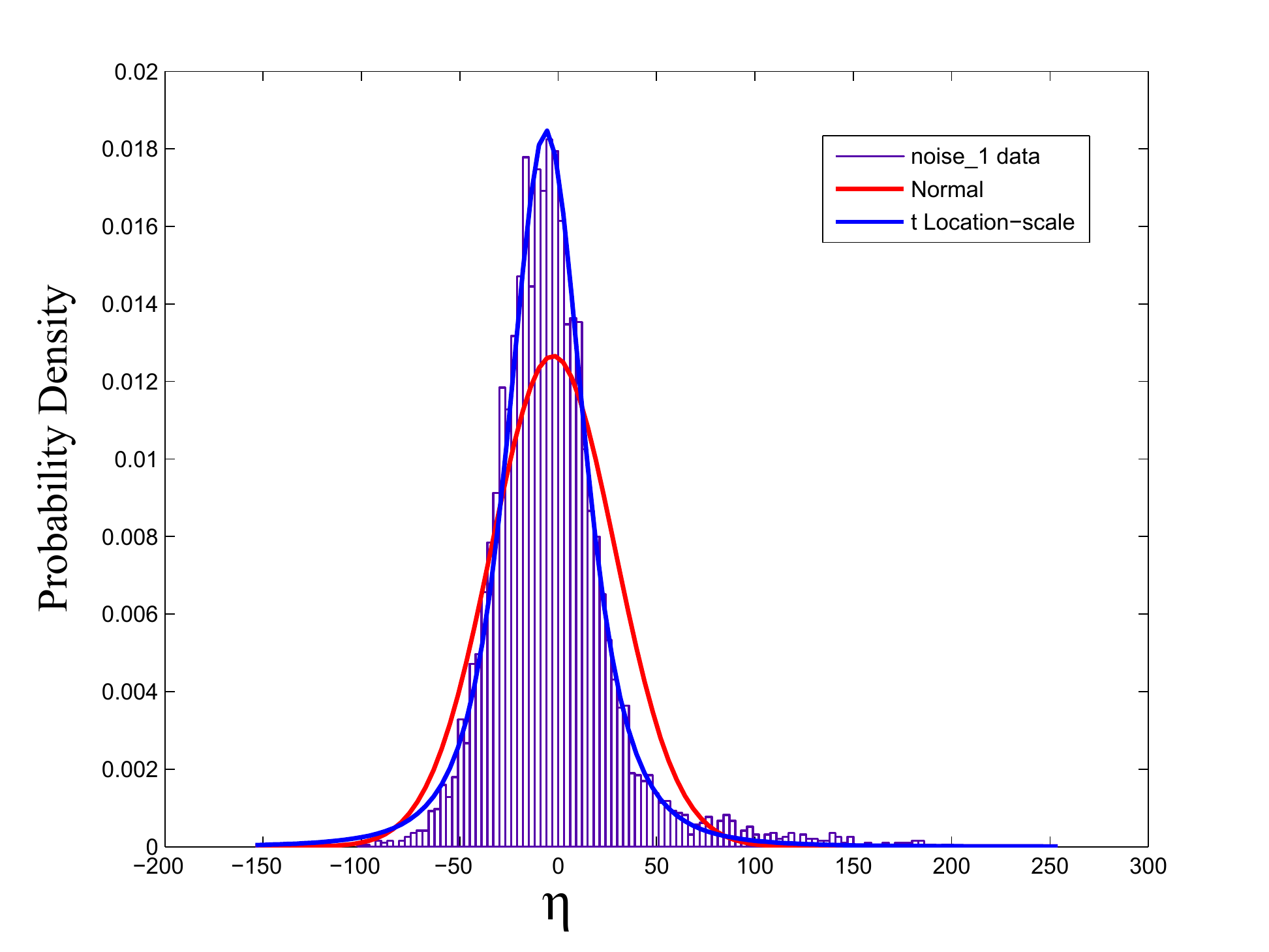}
\caption{Distribution of noise data points. A histogram of the data points in the time-series for the noise $\eta(t)$ 
(obtained using the cut-off $A_{\textrm{thresh}} = 20,000$)
is displayed along with fitted normal and $t$ Location-scale distributions. The mean and standard
deviation of the best-fit normal pdf are $\mu = -2.7$ and $\sigma = 32$. The parameters
of the best-fit $t$ Location-scale pdf are $(\mu,\sigma,\nu) = (-6.2,20,3)$.}
\label{fig:5}
\end{figure}
We now examine the noise $\eta(t)$, after removing the periodic signal using $A_{\textrm{thresh}} = 20,000$.
We assume that $\eta$ is a time-homogeneous sequence of random variables with a different realization for every $t$:
\begin{equation}
\textrm{Prob}[ x \leq \eta(t) \leq x+ dx ] = f(x) dx.
\end{equation}
Attempts at fitting this data to various well-known probability distributions reveals two distributions that yield a good fit: normal and the $t$ location-scale distributions (see Fig \ref{fig:5}) with  $t$-location-scale distribution yielding a better overall fit.
For a normal distribution, the we obtained the best-fit $N(\mu,\sigma^2)$ with $\mu = -2.7$ and $\sigma = 32$.
For the $t$ location-scale distribution with probability density function \begin{eqnarray}
f(x) &=& \frac{\Gamma[(\nu+1)/2]}{\sigma \sqrt{\nu \pi} \Gamma(\nu/2)}
\left[ \frac{[(x-\mu)/\sigma]^2+\nu}{\nu}\right]^{-\frac{\nu+1}{2}},
\end{eqnarray}
we found that the best-fit parameters were $\mu = -6.2$, $\sigma = 20$ and $\nu = 3$.
\subsubsection{Viewership as a stochastic process}
\begin{figure}[h!]
\begin{center}
\includegraphics[width = 0.6\textwidth]{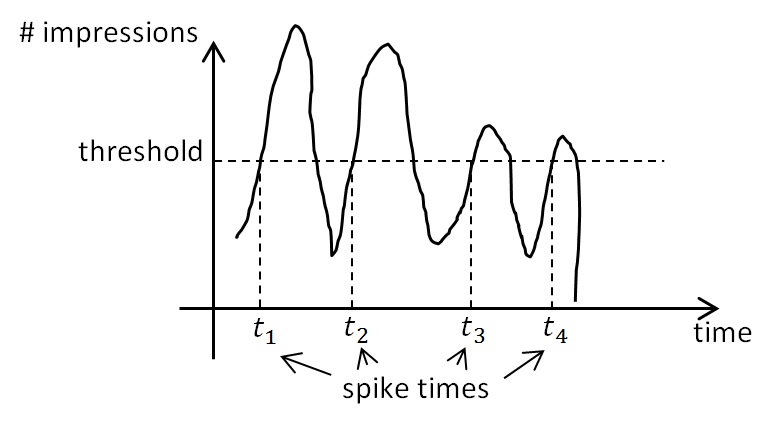}
\end{center}
\caption{Definition of spike time. As highlighted in the schematic, spike time is defined as the beginning of a series of impressions that exceed
a certain threshold. In our analysis, the threshold is set to be the 95th percentile of the data.}
\label{fig:sp1}
\end{figure}
\begin{figure}[!htb]
\centering 
\includegraphics[width = 0.45\textwidth]{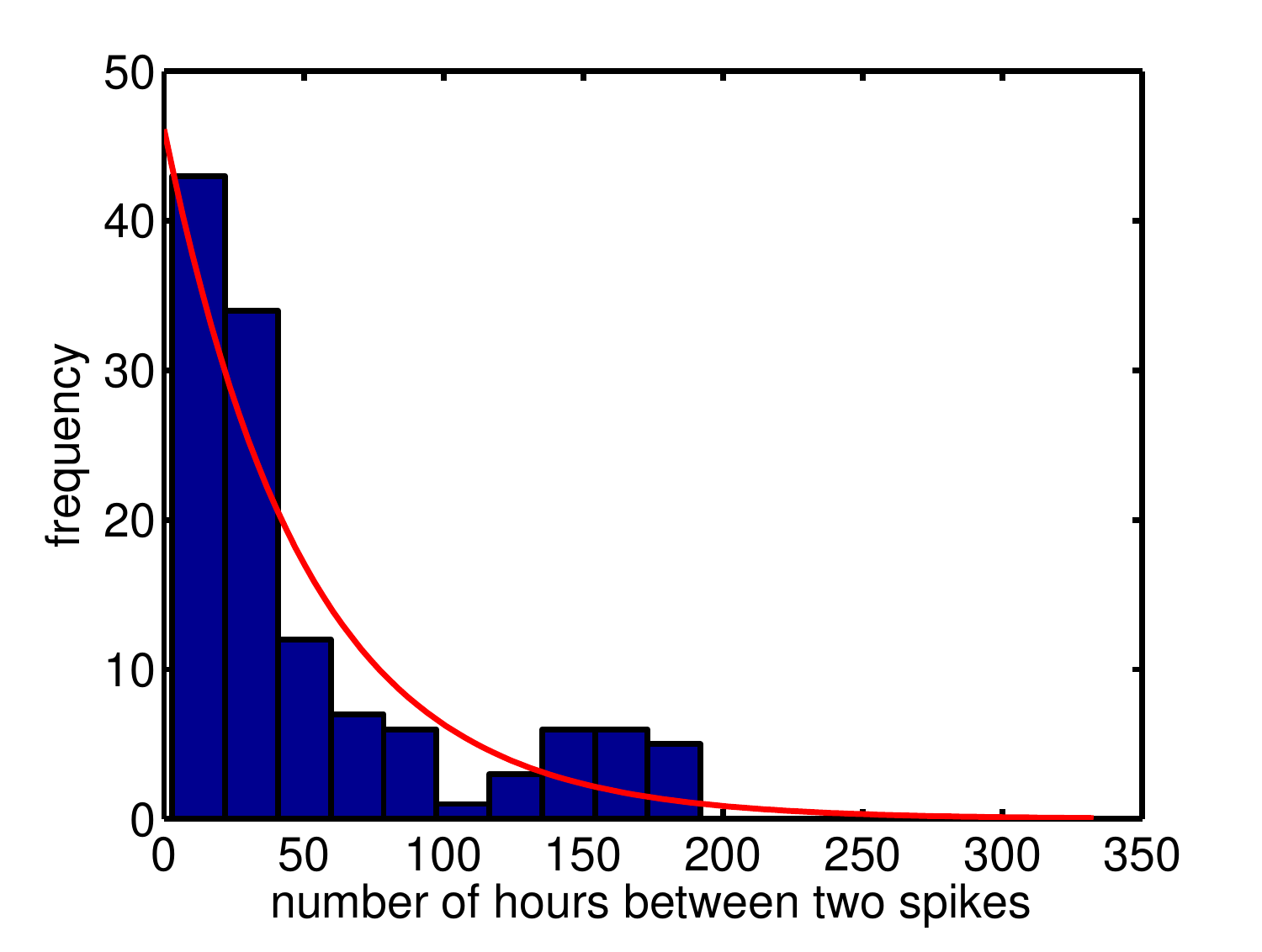}
\includegraphics[width = 0.45\textwidth]{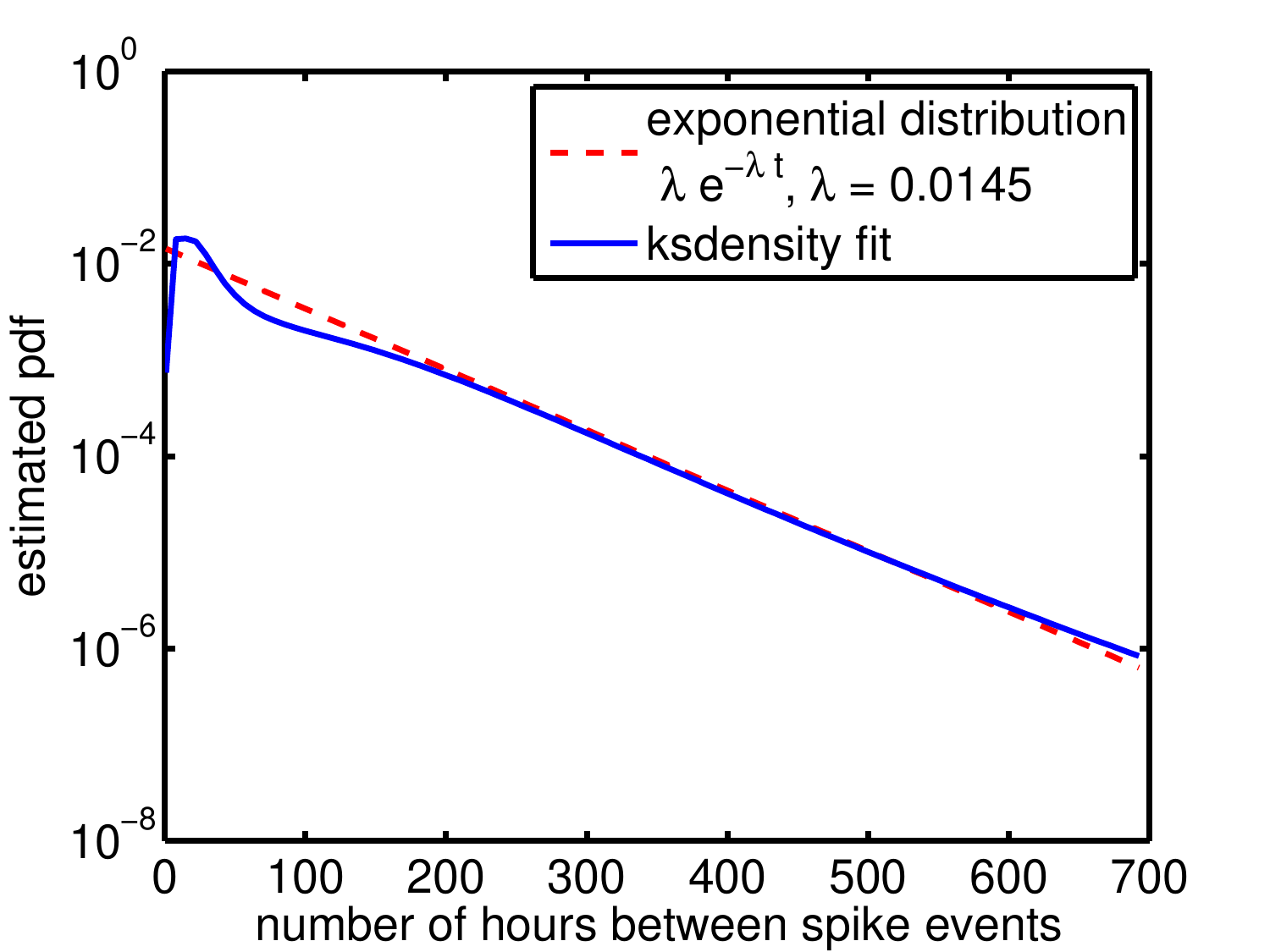}
\caption{Distribution of spike times. The time between two spikes in viewership approximately follows an exponential
distribution. \texttt{ksdensity} is a Matlab routine that attempts to find the pdf associated
with given data.}
\label{v4}
\end{figure}

One notable feature of $S(t)$ is that large spikes seem to occur randomly in time. We define the spiking event time as the time when the impressions signal
crosses a fixed threshold from below (see Fig. \ref{fig:sp1}).
Here we choose the threshold to be the 95th percentile of the available impressions data points. 
In Figure \ref{v4}, we show that the distribution of waiting times $\tau$ between consecutive spikes
appears to be approximately exponentially distributed:
\begin{equation}
\textrm{Prob}[ t \leq \tau \leq t+ dt ] = \lambda \exp(-\lambda t) dt.
\end{equation}
This suggests that the spiking has no memory
(spiking is approximately Markovian) and occurs at a Poisson rate of $\lambda \approx 0.015$
per hour, so that the mean time between spikes is about 69 hours. The analysis
of spiking time was performed on unfiltered data $S(t)$. One might also consider spikes in the noise left over from the filtering, $\eta(t)$. However, this yields similar results because the periodic signal generally has small amplitude.

\subsection{Machine Learning Approach}
Although spectral analysis of viewership data provides insight into the mechanisms that contribute to the observed trends, this approach assumes periodic behavior and ignores programming information. These findings can be helpful for filling in missing data and for estimating viewership of new programming, but an approach that takes into account the programming could potentially explain both the periodic behavior and the noise.  We therefore implement a machine learning algorithm to predict the number of impressions in a time slot by learning from past data. The machine learning task is defined with the following attributes: program ID (nominal),  day of the week (nominal) and time of the day (numeric). The output class is the number of impressions.

We explored two different approaches. First, we treated the output class as numeric. A number of machine learning methods are suitable for this task, and we tested k-nearest neighbor, neural networks, linear regression, regression tree, and k-star. The best performing algorithm was 1-nearest neighbor. The root relative square error for this method was 61\% (100\% would correspond to the error in naively guessing the mean of all impressions).

Second, we divided the output class into 5 bins, and tested 4 algorithms: decision trees, random forest, naive Bayes, and random tree. Decision tree and random forest performed equally well - the error was 44\% (compared with 80\% from naively guessing the correct bin).

\begin{figure}[h!]
\begin{center}
\includegraphics[width = 0.6\textwidth,clip=true, trim=0cm 0cm 0cm 0.7cm]{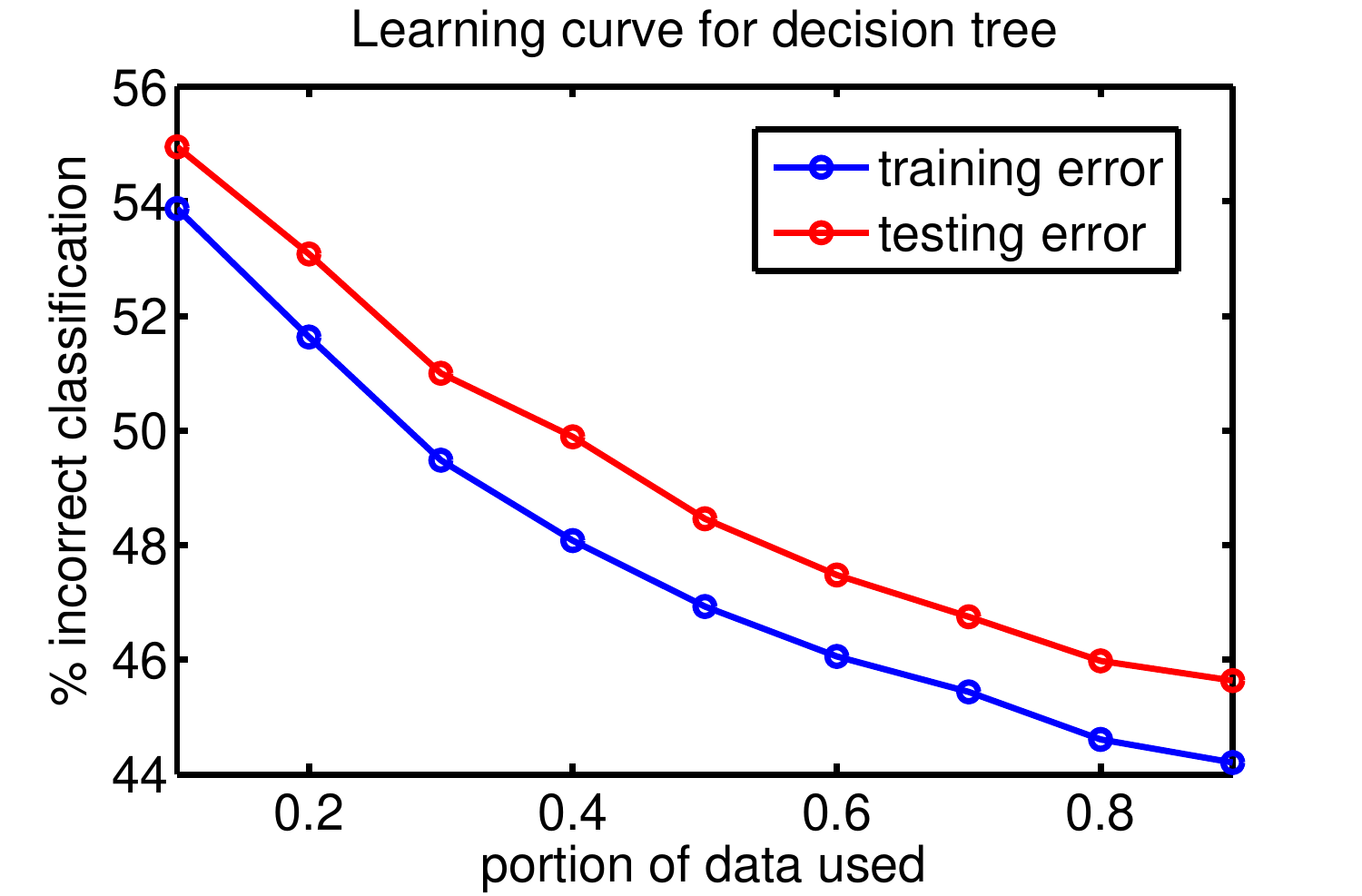}
\end{center}
\caption{Learning curve for decision tree machine learning algorithm. The lack of a plateau suggests that the algorithm could be improved with more data.  Given the small range of variability, alternative decision trees should be explored. }
\label{fig:ml1}
\end{figure}

A learning curve is shown in Figure~\ref{fig:ml1}. The fact that the curve did not plateau shows potential for more accurate prediction given more data. The fact that testing error decreases with training error demonstrates that this approach does not over-fit the data.

\subsection{Kalman Filtering}
\label{sec:kalman}
We also investigated the use of Bayesian estimation and Kalman filtering to predict the number of impressions for a program. In order to do this, we neglected any sampling issues that may be present in the  data and assumed perfect data. We treated the number of interested viewers as fixed with a particular probability of watching the program or not. This probability can be modeled by a binomial distribution because of the two possible values; however, due to the large sample size of viewers, we can apply the Central Limit theorem and approximate the distribution with a normal distribution \cite{Mood1974}. 
\par
For a Bayesian estimation, a prior probability distribution and likelihood function must be assumed in order to formulate an initial prediction. After this initial prediction is made, the new data is observed and is used to update the probability distribution and gives a posterior probability distribution. For our purposes, this posterior distribution was then used as our prior for predicting the next week \cite{Humpherys2012}.
\par
The number of impressions, $ S $, for each week $ w $ was modeled by a Gaussian with mean $ \mu (w) $ and standard deviation $ \sigma $. While we allow $ \mu $ to vary from week to week, we keep $ \sigma $ fixed for simplicity.   Under the Bayesian framework, we view $ \mu (w) $ as an unknown parameter which we will represent by a subjective probability distribution.  We can think of $\mu (w) $ as a measure of the popularity of the show at week $ w$, while the actual viewership will have an unpredictable fluctuation from week to week due to external factors.  The standard deviation $ \sigma $ measures this inherent variability in weekly viewership.
\par
To find a reasonable value for the fixed $\sigma$, the available data was divided into a particular number of bins. For each bin of data, the standard deviation was found, and then all of the standard deviations were averaged together to get an \emph{average} standard deviation. This was repeated multiple times with varying numbers of bins in order to choose the optimal (smallest) standard deviation. The smallest value found was then used as $\sigma$ for $S$.
\par
\begin{figure}[ht!]
\centering
\includegraphics[scale=0.62,clip=true,trim=1.5cm 2cm 1.5cm 1.62cm]{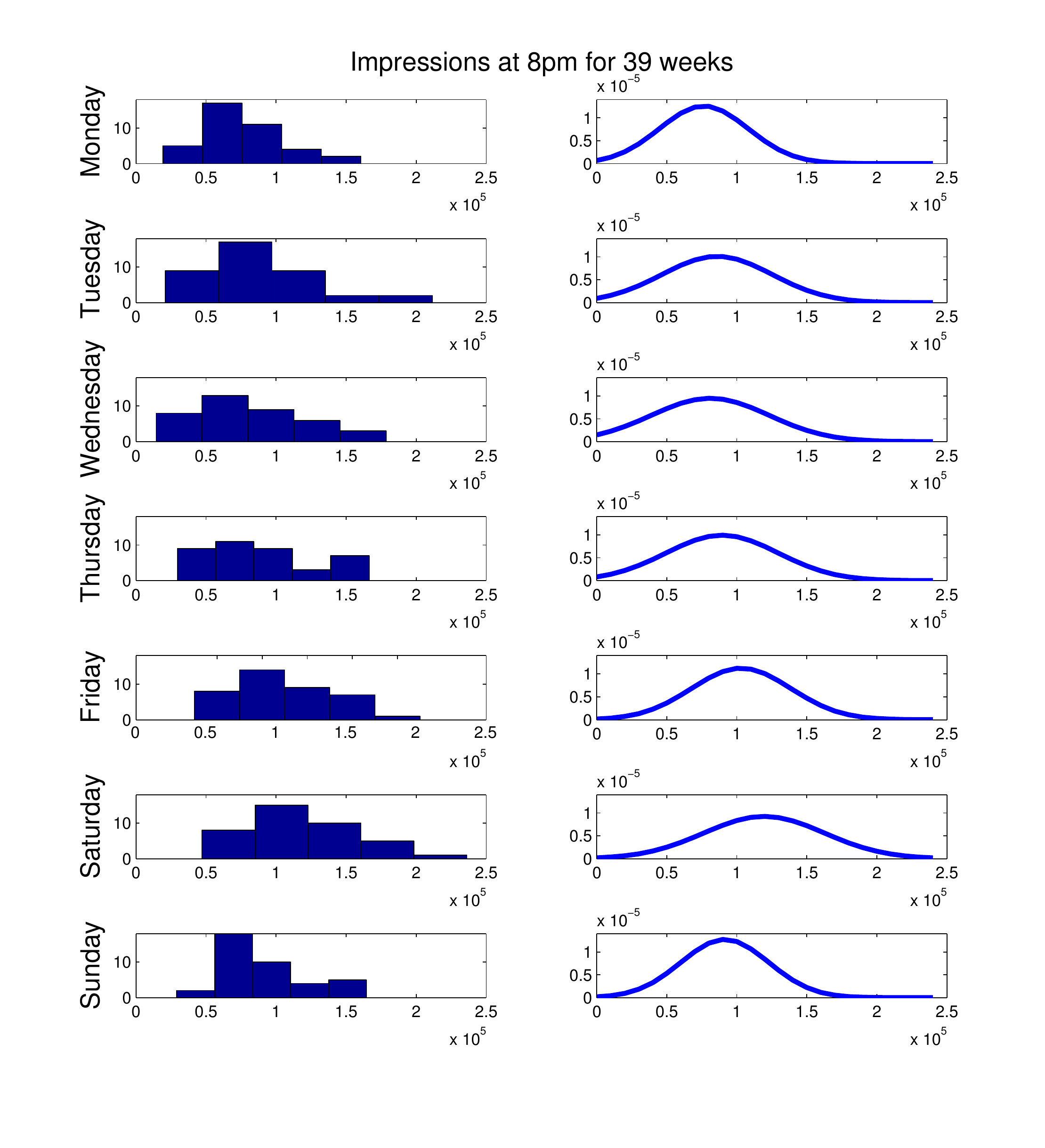}
\caption{Impressions by day. (Left) Histogram of the number of impressions during the 8 p.m. slot for each day of the week over 39 weeks. (Right) Approximating normal distribution of available data.}
\label{fig:histnorm}
\end{figure}

To understand the distribution of $\mu$, a recursive Bayesian estimation was used. To start, a weak prior probability distribution was chosen for $\mu$ based on a Gaussian fit to a histogram of 39 weeks of historical viewership data for a given time slot on a given channel on a given day.   These histograms and their Gaussian fits are illustrated in Figure~\ref{fig:histnorm} for 7 different weekly time slots, namely 8 p.m. on the given day of the week and a given channel.  In particular, the prior distribution for $ \mu $ will depend on the time of the week (and channel) under consideration.
Bayes' rule is applied in the usual way to update the prior distribution for an upcoming week with the likelihood of the data, once collected, to obtain a posterior distribution for $ \mu $ on that week.  This posterior distribution is then taken as the prior distribution for $ \mu $ on the following week.
The likelihood model, as described above, is Gaussian, so the recursive Bayesian estimation procedure reduces to a simple version of the Kalman filter, with the predicted (prior) distribution of the parameter $ \mu $ at the next week taken to be the same as the posterior distribution on the current week.  The evolving value of $ \mu $ is used as a point estimate for the expected number of impressions for a given channel on a given day of the week at a given time.
\par
In order to assess the performance of this model, we chose the first 20 weeks of data as the ``training'' data and then tested the model against the last 19 weeks. We tested 7 such data sets, namely the viewereship of a given channel at 8 p.m.~ on one of the 7 days of the week. For every day of the week the relative errors between the predicted impressions and observed number of impressions for the last 19 weeks were computed. The root mean square (RMS) of the relative errors was calculated as the measure of error for our model. The RMS error for each day of the week was below 30\%, where days such as Tuesday, Friday and Saturday were under 10\%. This suggests that our model is able to make reasonably accurate predictions for the number of future impressions given the data from the first 20 weeks. These results are displayed in Table \ref{impress}. The table also includes an example of our model's prediction for a week (39$^{th}$ week) for a particular network at 8:00pm.
\begin{table}[ht]
\centering
\begin{tabular}{|c|c|c|c|c|c|c|c|}
\hline
\textbf{Day} & \textbf{Mon}& \textbf{Tues}& \textbf{Wed}& \textbf{Thurs}& \textbf{Fri}& \textbf{Sat}& \textbf{Sun}\\[0.8ex]
\hline
\textbf{Date} & \textbf{2/2}& \textbf{2/3}& \textbf{2/4}& \textbf{2/5}& \textbf{2/6}& \textbf{2/7}& \textbf{2/8}\\[0.8ex]
\hline
\textbf{Model} & 83902& 101720& 65465& 50160& 65032& 101120& 63138\\[0.8ex]
\hline
\textbf{Actual} & 87204& 105940& 62146& 39267& 60376& 93997& 53142\\[0.8ex]
\hline
\textbf{Rel. Error} & 0.0379& 0.0398& 0.0534& 0.277& 0.0771& 0.0758 & 0.188\\[0.5ex]
\hline
\textbf{RMS for Day} & 0.27& 0.0903& 0.2203& 0.2016& 0.0804& 0.0883& 0.1983\\[0.5ex]
\hline
\end{tabular}\\
\begin{tabular}{c}
\\
\end{tabular}
\caption{Impressions estimate for February 2-8, 2015 for a given network at 8:00pm. The model estimates are given along with the observed number of impressions. The relative errors for the displayed week is calculated and compared to the root mean square determined from all of the 39 weeks.}
\label{impress}
\end{table}
\subsection{Distance scores}
One shortcoming the previous two approaches is that they require accurate data about the viewers of a particular program to train the algorithm.  However, for new programming it is necessary to generate predictions for the number of impressions with limited data. To overcome this challenge, it is helpful to consider the viewership trends of similar programs.  We therefore created a difference score to identify appropriate similar programs for comparison. Given two programs, $p_1$ and $p_2$, at times $t_1$ and $t_2$ respectively, we define a function $\Delta(p_1,t_1,p_2,t_2)$ that returns a score measuring how different those programs are at those times in terms of their demographic ratios. The idea is that the viewership breakdown of future television shows can be predicted by studying the breakdown of similar shows that have aired in the past.

First we let $D_{p,t}$ be the demographics information for program $p$ at time $t$ (where $t$ includes information about day of the week and time slot during the day, ex: Monday from 3pm to 4pm).
$$D_{p,t} = (d_1,d_2,...,d_N)$$
where $d_i$ is the number of impressions from the $i$-th demographic (ex: females ages 21 to 24) and $N$ is the number of demographic fields in the data (in our case $N=30$, with $15$ age ranges for both males and females). Now let $I_{p,t} =\sum_{i=1}^N d_i$, be the total number of impressions for program $p$ in the time slot $t$. Let $\hat{D}_{p,t}=\frac{D_{p,t}}{||D_{p,t}||_1} = D_{p,t}/I_{p,t}$ which is $D_{p,t}$ normalized in the $\ell_1$ norm. $\hat{D}_{p,t} \in \mathbb{R}^N$ is a  breakdown of the viewership by demographic for program $p$ at time $t$.

The distance score is defined as 
\begin{equation}
\Delta (p_1,t_1,p_2,t_2) = ||\hat{D}_{p_1,t_1} - \hat{D}_{p_2,t_2}||_2,
\label{eqn:S}
\end{equation}
which is just the Euclidean distance between $\hat{D}_{p_1,t_1}$ and $\hat{D}_{p_2,t_2}$ in $\mathbb{R}^N$. Note that this means that programs with lower scores are more similar because their viewer demographics will be closer to each other in terms of the Euclidean distance. 

Since we do not know where in the space our data sits, we randomly selected 
pairs of programs to see how far apart they are on average, keeping certain program attributes (such as
program ID, day of showing, hour of showing) identical. Our results are shown
in Table \ref{tab:sim2}. This reveals that the programming content is a better predictor of the demographic data than the programming date and time.

This difference measure relies on the $L^2$ norm, however alternative metrics might also yield reasonable scores. For example, we could compute the difference between $D_{p,t}$ and $D_{p', t'}$ using a dot product
\[ \langle \hat{D}_{p,t}, \hat{D}_{p',t'} \rangle\]
resulting in a score between $0$ and $1$ (after dividing by a normalization factor).

These scores could be used to estimate viewership of new programs by averaging the impression data from a group of similar programs. These predictions could then be enhanced by correcting for the spectral properties of the viewing habits of each demographic in the relevant timeslots. Then as additional data is collected, these predictions could be adjusted using the Kalman filtering approach.

%

\begin{table}
\centering
\begin{tabular}{l l l c}
\hline
\textbf{Program ID} & \textbf{Day of showing} & \textbf{Hour of showing} & $
\Delta$\\
\hline
Random & Random & Random & 0.5128\\
Random & Same & Same & 0.4845\\
Same & Random & Same  & 0.3146  \\
Same  & Same  & Random  &  0.2842 \\
\hline
\end{tabular}\\
\begin{tabular}{c}
\\
\end{tabular}
\caption{Average distance score for 100,000 randomly selected program pairings. 
The distance $\Delta$ is defined in equation (\ref{eqn:S}).}
\label{tab:sim2}
\end{table}

\section{Optimal Scheduling using Integer Programming}
\label{sec:BIP}
In this section, we implement a method for creating an optimal schedule of advertisements, given
a set of orders from an advertising agency and predicted viewership numbers such as those that could be generated using the methods outlined in Section \ref{sec:predict}. 


In order to formalize the optimization method we define the following notation and assumptions:
\begin{itemize}
\item $C$ is the total number of channels, and the subscript index $c$ with $1\leq c\leq C$ is used to denote one particular channel.
\item $N_c$ is the number of commercial slots on channel $c$, and the subscript index $i$ with $1\leq i\leq N_c$ is used to denote the corresponding slot. This index takes into account both day and time.  We assume the number of slots are specified by programmers in advance.
\item $\Pci$ indicates the price for commercial slot $i$ on channel $c$. We assume these prices are set by the programmer in advance.
\item $\Scid$ contains the number of impressions for slot $i$, demographic group $d$, on channel $c$. We assume these values are provided in advance and that this data is reliable.
\item $A$ gives the number of advertising orders, and the superscript index $(a)$ with $1\leq a\leq A$ denotes one particular order. We assume that all orders for a given week are received in advance, that the schedule can be determined one week at a time, and that all advertisers have equality priority and therefore orders accepted or rejected only on the basis of whether the order is likely to be satisfiable.
\item $\Va$ is a binary vector indicating the target demographics in the order for advertiser $a$. 
\item $\Scia$ contains the number of impressions for slot $i$, in the demographics specified by advertiser $a$, on channel $c$. In other words, $S_{c,i}^{(a)} = \sum_{d \in \mathbf{V}^{(a)}} S_{c,i}^{(d)}$. We assume that all target demographics are of equal value to the advertiser and therefore the desired number of impressions can be satisfied by any subset of the target audience.
\item $\Ba$ represents the budget of advertising order $a$, and $\Ra$ represents desired impressions for order $a$. We assume that these requirements are strict and can be implemented as hard inequality constraints for the solution.
\item $\Xcia$ is a `binary matrix' indicating whether advertiser $a$ is assigned to slot $i$ on channel $c$. This is the schedule we are trying to find.
\end{itemize}
\subsection{Problem Formulation}
We now use this notation to express the scheduling problem as a constrained optimization problem.
\subsubsection{Constraints}
The most basic constraint on a proposed schedule is that two advertisements cannot air simultaneously on the same channel.\\
\noindent \textit{1.  No overlap:}\\ Only one advertiser can use a given slot on a given channel. Mathematically, this can be stated as
\begin{equation*}
\sum_{a}\Xcia\leq 1.
\end{equation*}
This constraint can be modified to allow for variable length commercials by weighting each entry in $\Xcia$ by the commercial length for advertiser $a$ and then changing the right hand side to include the number of `time slots' in each commercial break.

In addition to this, each order $(a)$ that is accepted imposes two additional inequality constraints on the schedule.  \\
\noindent \textit{2.  Budget:}\\ The total cost to each advertiser must not exceed their budget $\Ba$. This implies that
\begin{equation*}
\sum_{c,i}\Xcia\Pci\leq \Ba.
\end{equation*}
Note that it may not be possible to satisfy every order, so if the total cost to an advertiser is greater than 0, then we must also meet the target number of impressions.\\
\noindent \textit{3. Impression Target:}\\ The total number of impressions (as given by $\Scia$) must exceed the campaign goal $\Ra$.  In other words,  
\begin{equation*}
\sum_{c,i} \Xcia S_{c,i}^{(a)} \geq \Ra.
\end{equation*}
Since this linear inequality only yields a feasible region if it is possible to satisfy every order (which in practice is unlikely), we impose these constraints by solving a sequence of optimization problems where $\Ra$ are replaced with 0 for the orders we are not able to fill. In Section \ref{vf}, we propose a value function that can be used to determine which orders should be eliminated. 

The above constraints are necessary to obtain a usable schedule that satisfies the advertising campaign goals. However, programmers may impose additional requirements on allowable schedules to prevent consecutive airings of the same commercial ($\Xcia+X_{c,i+1}^{(a)}\leq 1$ for all $i$, $c$, $a$), commercials with adult content from airing during children's programming ($\Xcia=0$ for particular $i$, $c$, $a$), etc.  These and any other requirements can be implemented by imposing additional inequality and equality constraints on $\Xcia$.  However, for simplicity, we omit these constraints in what follows.
\subsubsection{Objective function}
Presumably, the advertising schedule is set by the programmer or an intermediary who is interested in maximizing advertising revenue.  Therefore, the objective function of interest is simply the total revenue which is given by 
\begin{equation}
\sum_a\sum_{c,i}\Xcia\Pci.
\end{equation}
\subsubsection{Binary Integer Program}
This optimization problem involves finding a vector inputs $\mathbf{X}$ (a vectorized version of $\Xcia$) that satisfies a set of linear constraints and that maximizes the value of a linear objective function. If the inputs were real numbers, then this could be solved with a linear program.  However, $\mathbf{X}$ must contain binary inputs so therefore we solved this using a binary integer program. 

To implement this program, we write the inequality constraints as a matrix inequality \[A\mathbf{X}\leq B\] and the objective function as dot product \[f(\mathbf{X})=\mathbf{P}\cdot\mathbf{X},\] 
where $\mathbf{P}$ is a vectorized version of $P_{c,i}$. This allows us to make use of MATLAB's built-in mixed integer linear programming algorithm from the optimization toolbox. This algorithm consists of the following three steps \cite{Wolsey1998,Wolsey2014,Danna2005}:
\begin{enumerate}[1.]
\item Solve the linear programming problem without the integer valued constraints.
\item Use a heuristic algorithm to find a nearby feasible integer solution.
\item Perform branching to try to improve on the heurtistic feasible solution.
\end{enumerate}
Depending on the data and parameters used, this algorithm occasionally finds a solution which sells all of the time slots or it fills all of the order leaving some time slots unfilled. These outcomes represent the global maximum for the revenue. Other times, it cannot find a feasible solution at all.  This suggests that a high quality heuristic is necessary for finding feasible solutions to initialize the integer program \cite{Savelsbergh1994}. We explore one such heuristic in Section \ref{greedy}.  One explanation for this inability to find a feasible solution is the fact that it is not always possible to satisfy all of the orders.  To overcome this, we iteratively remove orders and instead choose a subset of the orders that includes only the ones that are the most valuable (see Section \ref{vf}). 

\subsection{Greedy Algorithm}\label{greedy}
In order to generate a feasible solution both to initialize the integer program and to compare to the results from the integer program, we also implemented a greedy algorithm.  In this algorithm, we generate a matrix $V$, where the rows are indexed by the slot $\{i,c\}$  and the columns are index by the advertiser $a$. We assign a value to each entry of $V$ for each advertiser.  Then we choose the slot with the highest value and assign that to the advertiser who gets the most value from that slot.  The value of slot $\{i,c\}$ is for advertiser $a$ is equal to
\begin{equation}
V_{c,i}^{(a)}=\left(\frac{\Scia}{\Ra}\right)/\left(\frac{\Pci}{\Ba}\right)
\end{equation}
which represents the fraction of the desired impressions that can be provided by the slot divided by the fraction of the budget that must be used for the slot. 

This process is repeated until there are no longer available slots, the orders have all been met or no advertisers can afford a slot. At this point, incomplete orders are removed and the process is repeated with the remaining orders until all orders have been satisfied or removed. 
\subsection{Value function for individual orders}\label{vf}
In order to determine which orders should be rejected, a systematic way of prioritizing orders is needed. Below, we propose a heuristic ranking scheme.\\\\
\noindent \textit{Step 1: Eliminate unreasonable orders}
\\\noindent If the number of impressions desired is more than the size of the viewership for that demographic (in the time allotted), then the order cannot be satisfied. These orders should be rejected immediately.  In other words, an order must be rejected if  $\sum_{c,i}\Scia<\Ra$.\\\\
\noindent \textit{Step 2: Monte Carlo Method}
\\\noindent A Monte Carlo method can be used to estimate the number of feasible solutions for each advertiser. We then assign a value proportional to that number.
\begin{enumerate}[1.]
\item First generate a random binary vector $\Xcia$ for fixed $a$.
\item Check to see if it is feasible.
\item Repeat $N$ times.
\end{enumerate}
Let $\mathbf{F}$ be an $N$ by 1 vector indicating which candidate solutions are feasible.  Then, the fraction of solutions that are feasible is given by
\begin{equation}
W^{(a)}_{MC}=\frac{\mathbf{F}\cdot\mathbf{1}}{N}
\end{equation}
which we call the value function, and the orders that are more likely to be satisfiable should be prioritized.
\subsection{Extensions of the value function}
This provides a simple method for identifying feasible orders. However, in practice the decision making process might be more complex.  For example, if rather than having fixed time-slot prices advertisers were allowed to bid for a given slot, then the advertisers with larger budgets relative to their demands should be prioritized since their orders are likely to be both more satisfiable and more lucrative.  Also, if a variable scheduling horizon for orders is allowed (rather than focusing on a single week at a time), then the urgency of the order should also be taken into account. \\\\ 
\noindent \textit{Modification 1: Weight by excess budget.}\\\noindent 
This value function can be modified to account for bidding on time slots by weighting the feasible solutions by the price the advertiser would
be willing to pay in an auction.  An advertiser should be willing to increase the bid by a factor of $\frac{\Ba}{\sum_{c,i}\Xcia\Pci}$ to remain under budget and ensure that their advertising campaign order is accepted. So, given a set of random $\Xcia$ from the Monte Carlo method above, rather than just computing the fraction that are feasible, the feasible schedules can be weighted to account for the amount of leftover money in the budget. For a feasible solution $j$, the percentage of the budget that is unused is just 
\begin{equation*}
E^{(a)}_j=\frac{\Ba-\sum_{c,i}\Xcia\Pci}{\Ba}.
\end{equation*}
Let $\mathbf{F}$ be an $N$ by 1 vector indicating which candidate solutions are feasible.  Let $\mathbf{E}^{(a)}$ be the percentages from above. Then the total value of an order would be given by the average of this excess 
\begin{equation*}
W^{(a)}_{BID}=\frac{\mathbf{F}\cdot\mathbf{E}^{(a)}}{N}
\end{equation*}
and orders with larger average value should be prioritized.  In order to balance both feasibility and value, a balance of the Monte Carlo and bidding based values such as  
\begin{equation*}
W^{(a)}_{COMB}=(1-r)W^{(a)}_{MC}+rW^{(a)}_{BID}
\end{equation*}
was employed, where $r$ represents the relative weight of the excess budget to the feasibility. \\\\
\noindent \textit{Modification 2: Weight by urgency.}\\
In order to take into account variable time frames for orders, this value could
factor in the fact that certain orders are more urgent than others. Orders that have already been accepted will usually need to be prioritized over new orders. For new orders, some orders with short time tables will be infeasible, others with short time tables would need to be completed immediately, and orders with long time tables may be saved for later, but not postponed so long that they become infeasible. This modification was not implemented, but it warrants further exploration.
\subsection{Results}
In order to test this optimization scheme, we used a subset of the viewership data over a single work week (Monday-Friday) between the hours 5:00am-12:00am for three channels. We also assumed that impressions by demographic are constant over the span of an hour (if there are two half hour shows in an hour we average the impressions per demographic from both shows) and do not take into account any uncertainty in the estimates of impressions per time slot.  The results are displayed below.
\begin{figure}[h!]
\begin{centering}
\includegraphics[width=.65\textwidth]{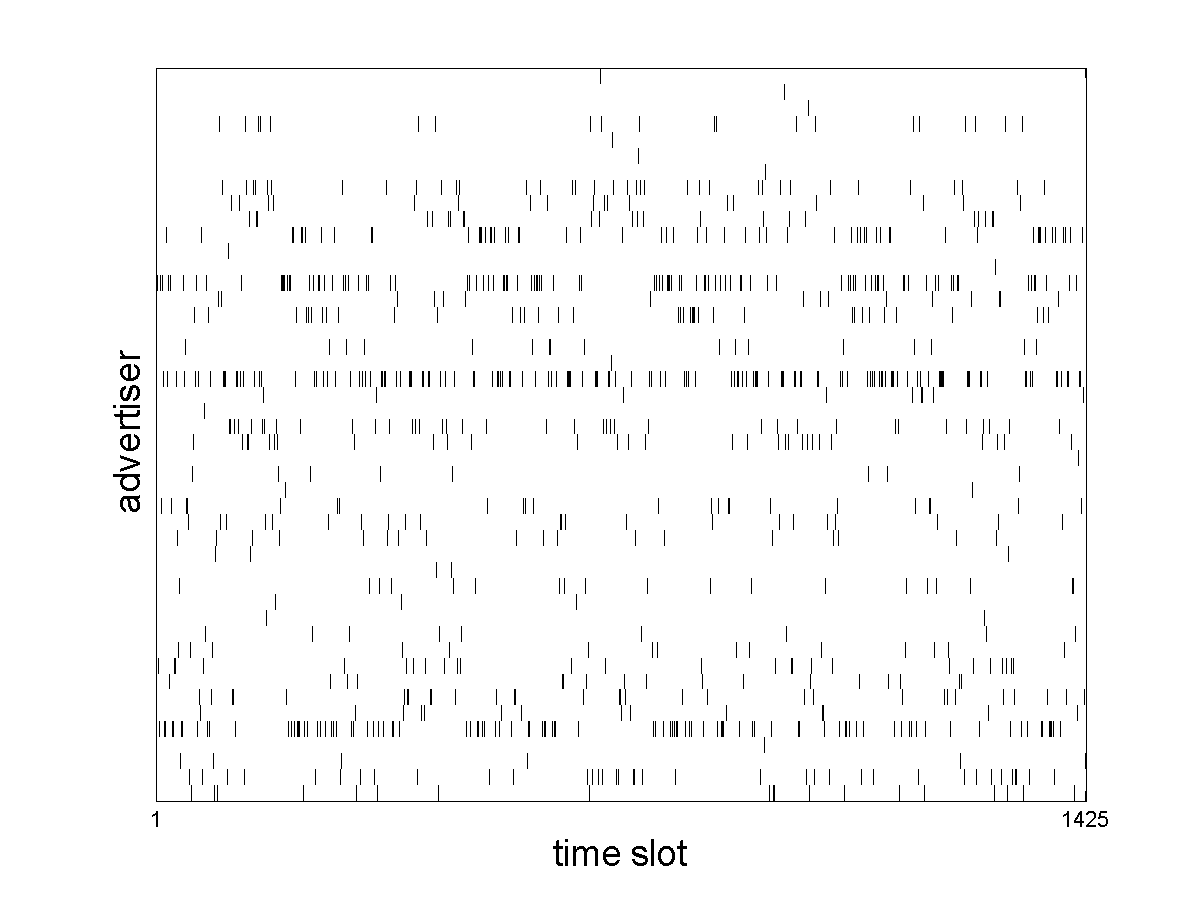}
\caption{Optimal schedule schedule for a programmer with one slot per hour on three different channels for an entire work week. The resulting schedule satifies the requirements for the top 49 orders over this span. The vertical axis corresponds to different advertiser's orders while the horizontal axis corresponds to the time slot index. Black corresponds to a slot being filled.}
\label{sample_schedule:fig}
\end{centering}
\end{figure}

For this data set, an optimal schedule can be found that satisfies the top 49 orders and fills all available ad time when orders are sorted using the Monte Carlo value function. The optimal schedule is shown in Figure \ref{sample_schedule:fig}. With more orders, however, the algorithm fails because it is unable to find a feasible solution. Thus iterative reductions in the total number of orders are necessary before a satisfiable subset of orders can be found.  In contrast, the greedy algorithm always yields a feasible solution (by automatically rejecting unsatisfied orders), but that solution may not be close to optimal.  With 149 orders (the total number of sample orders in our data-set), a solution satisfying 112 orders that generates 91.8\% of the maximum allowable revenue is found in 0.57 seconds whereas the integer program finds the optimal solution satisfying as many as 49 orders in 6.5 seconds.

\begin{figure}[h!]
\begin{centering}
\includegraphics[width=.65\textwidth]{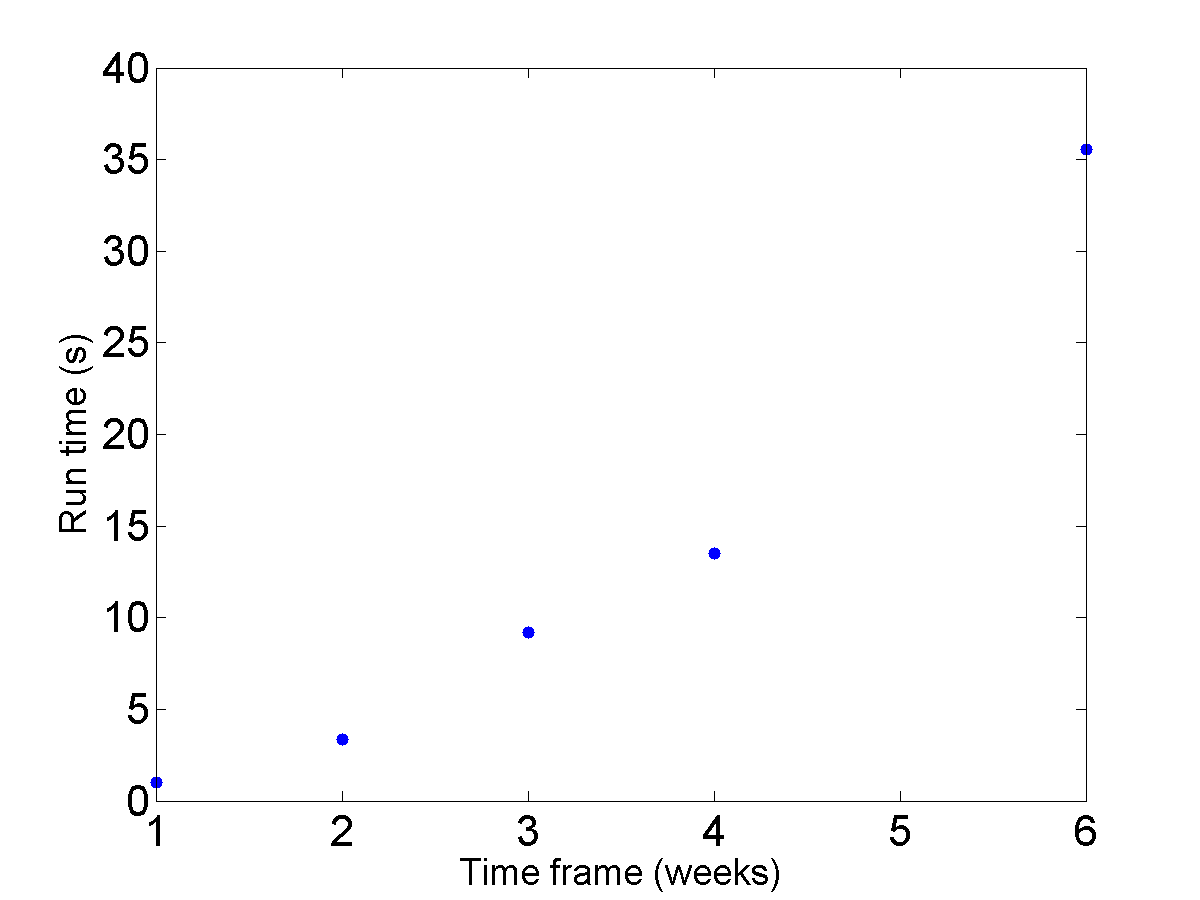}
\caption{Run-time for computing an optimal schedule with a variable advertising time frame. The integer program was initialized with a variable number of weeks and fixed set of 40 orders (rescaled to match the number of weeks). These data points can be fit by a quadratic $y=0.96535x^2$ with $R^2=0.994$.  This suggests that the computation time is proportional to the square of the time frame.}
\label{run_time:fig}
\end{centering}
\end{figure}

In the our data set, adding more time slots by expanding the time horizon to several weeks or using more channels with varied viewership allows us to accommodate more orders, but it can also increase the computation time (see Figure \ref{run_time:fig}). Thus this method may be more suitable for smaller problems with fixed time horizons. Scaling the algorithm to larger data sets and to include more complex constraints (for example, to take into account the reach and frequency of an advertisement) would require a hybrid approach involving multiple algorithmic frameworks. When long time horizons, large numbers of channels or large numbers of orders must be considered, one promising approach would be to segment orders and schedules into smaller intervals and then apply this integer programming method to each interval.  The most efficient method for this segmentation would likely be dependent on the data considered but there may be structural properties of this type of problem that can be exploited. Therefore, scaling this method effectively would require a deeper look at segmentation strategies that would allow for the large scale problem to be divided into pieces that could be solved in parallel. 

These results suggest that binary integer programming provides a flexible framework for implementing the essential constraints and searching for optimal solutions. However, the development of new heuristics and improvement of current heuristics could be beneficial for scaling of the problem to more realistic sizes and warrants further exploration.

\section{Reach and Frequency}
\label{sec:reach}
In addition to the total number of impressions, advertisers may also be interested in specifying the desired \emph{reach} and \emph{frequency} of their advertisements.  The \emph{reach} of an advertisement is the number of unique individuals who have received at least one impression of the advertisement, and the \emph{frequency} of an advertisement is the mean number of times the advertisement is seen by these individuals.  From detailed viewership data, we can express the reach of an advertisement by the following exact formula:
\begin{equation}
R = \sum_{j=1}^{N^{(a)}} S^{\sharp}_{i_j} \label{eq:reach_exact}
\end{equation}
where $ S^\sharp_{i_j}  $ is number of \emph{new} impressions made at
time slot $ i_j $ of the $j$th airing of advertisement, and $ N^{(a)} $ is the total number of times the advertisement is shown.  To reduce the notational complexity, we are dropping here the indices referring to the channel and demographic group, essentially assuming we are focusing on a fixed channel and a specified demographic group.  This can be generalized in principle to allow multiple channels and multiple demographic groups, with accompanying complication in notation. 

The formula for reach is to be contrasted with the formula for total number of impressions made by the advertising campaign,
\begin{equation*}
I = \sum_{j=1}^{N^{(a)}} S_{i_j}
\end{equation*}
where $ S_{i_j} $ is the number of impressions (including both new and repeat viewers) made at the time slot $ i_j $ of the $j$th airing of the advertisement.  $ S^{\sharp}_{i_j} $ by contrast only counts those viewers on whom an impression was made at time slot $ i_j $, but not on a previous airing of that advertisement.   Consequently, to count reach, we must know more about the viewership than simply the statistics for the number of impressions likely to be made in each time slot.  Put another way, the number of impressions $ S_i $ is only a function of the data at time slot $ i $ (and thus may be thought of as a one-dimensional marginal distribution of the viewership data), whereas the number of new impressions $ S^\sharp_i $ depends not only on the statistics of time slot $ i $ by also on statistics of previous time slots (and thus inherently involves joint distributions of viewership data at different time slots).  To make matters more complicated for the purpose of schedule optimization, the number of impressions $ S_i $ in a time slot depends only on the time slot in which the ad is scheduled, whereas the number of new impressions $ S^{\sharp}_i $ depends not only on the time slot $ i $ but also the previously scheduled time slots of the advertisement.  The average frequency fortunately is easily determined from the reach by the simple formula:
\begin{equation*}
F = I/R.
\end{equation*}

\subsection{Predictive Scheme for Reach}
\label{sec:reach_predict}
As noted above, the reach of a previously aired advertisement can easily be calculated from historical data. However, predicting the reach of a proposed future advertising campaign is a much more delicate matter.  Even if the viewership of future programs could be assumed to be identical to previous weeks, the reach could be calculated for any proposed schedule, but this would be an expensive and unwieldy calculation. Either it would be necessary to do an online processing of historical data, or it would be necessary to refer to an intractably large data structure which has precomputed reach scores for every feasible advertising schedule.   Including reach into the optimization scheme for the advertising schedule would, as we show, introduce an inherent nonlinearity, and nonlinear optimization typically requires additional iterations beyond that required for linear optimization. That means many schedules will be proposed in the optimization, and thus the expensive reach computation would be invoked many times.  We therefore consider a simplified way to estimate future reach.  Such simplification is justified in particular because future viewership cannot be perfectly predicted, so involving a expensive and precise computation for reach in the schedule optimization algorithm would seem to be a misplaced effort.

We propose encoding the information needed for future reach calculations in a two-slot function $ P_{i,i^\prime} $ which, 
 for $ i < i^\prime $, represents an estimate of the fraction of  viewers  at time slot $ i^\prime $ who also viewed time slot $ i $.  Under our standing simpilfying assumptions, $ P_{i,i^\prime} $ could be estimated from historical data, possibly using the Kalman filtering idea in Section \ref{sec:kalman}.  If we allow future time slots to be associated with programs different from those in the past, we could try to develop an inference scheme for combining historical data on viewership of time slots with viewership of programs.  But this might still be attempting too fine a resolution.  Since the number of potential time slots in which a proposed advertisement could air within a typical campaign window is large, such a detailed data-driven approach would require the storage of $P $ as an immense matrix, which would be at least $ 10^3 \times 10^3 $ for a weeklong campaign even in our very simplified setting, and much larger in practice.  A more tractable approach might be to simply treat $ P_{i,i^\prime} $ as a function of only $ i-i^\prime $, meaning essentially the time difference between slots (and possibly also a measure of difference between channels for a multichannel campaign).  We might imagine that $ P_{i,i\prime} $ begins as a decaying function of $|i-i^\prime|$, but has peaks at multiples of a day and a week for patterned viewer behavior.   Perhaps historical data could be fit to a sum of a small number of periodic functions with frequencies identified by the spectral analysis in Section \ref{sec:spectral}, with decaying amplitudes.  

We now assume we have in hand some scheme for estimating the two-slot function $ P_{i,i^\prime} $, and now wish to estimate the reach of a proposed scheduling of the advertising campaign in time slots $ \{i_1,i_2,\ldots,i_{N^{(a)}}\} $.  According to formula~\eqref{eq:reach_exact}, we need to estimate the number of new impressions made with each airing of the advertisement, and we propose to approximate this in terms of the estimated number of impressions and the two-slot function as follows:
\begin{equation}
S^\sharp_{i_j} \approx S_{i_j} \prod_{j^\prime<j} (1 -P_{i_{j^\prime},i_j}). \label{eq:newimp}
\end{equation}
 That is, the estimated number of new impressions is equal to the estimated number of impressions, discounted by factors $(1 -P_{i_{j^\prime},i_j}) $ representing the fraction of the viewers of the $j$th airing of the ad who did \emph{not} also see the prior $j^\prime$th airing of the ad.  The approximation in Eq.~\eqref{eq:newimp} is conditional independence of the viewership of all previous airings of the advertisement by the viewers of the $j$th airing of the advertisement.  For a concrete example, for $ j=3 $, the \emph{conditional} independence assumption is that whether viewers of the third airing of the advertisement watched the first airing of the advertisement is independent of whether they watched the second airing.  Note that this is not the same as stating that a general viewer has an independent chance of viewing the first and second airing of the advertisement (unconditional independence).  Indeed, if $ P_{i_1,i_3} $ and $ P_{i_2,i_3} $ were $0.95 $, then the probability model underlying the formula~\eqref{eq:newimp} would have a substantial positive correlation between the viewers of the first and second airing of the advertisement.  The point of the conditional independence assumption is that we assume all such correlations between the viewership of the various airings of the advertisement can be well represented by an explicit model of the correlation between the viewer of each airing $ j^\prime < j $ and the airing $ j $ under consideration, with the correlations between the previous airings being implied (not neglected) by the conditional independence assumption.  

The conditional independence assumption can lead to either overestimates or underestimates of the reach.  For example, if the airings occur during successive episodes of a program with a substantial committed core base who watches every episode, the number of new impressions would be underestimated by formula~\eqref{eq:newimp}.  On the other hand, if the airings of the advertisement involve some episodes repeated at different times during a week, the viewership of those airings would be more negatively correlated than the conditional independence assumption, and the number of new impressions could be overestimated by formula~\eqref{eq:newimp}.   This can be verified under a simple model in which no viewer makes repeated viewings of the same episode at different times.  Some kind of conditional independence assumption seems to be necessary to reduce the reach calculation to a complexity comparable to the 2-slot statistic.  Another natural way to invoke conditional independence is via a Markov chain model, which would only attempt to explicitly model the repetition in viewership between successive airings of the advertisement.  Such a Markovian approach appears less suitable than the conditional independence we suggest in the previous paragraph for a couple of reasons.  First of all, it is unclear how to deduce the number of new impressions made on the third airing of an advertisement by knowing how many viewers of the first advertisement saw the second advertisement, and how many viewers of the second advertisement saw the third advertisement.  How does one infer from this the number of viewers of the third advertisement who saw neither the first nor the second airing?  Moreover, the Markovian approach seems completely incapable of representing the likely strong repeat viewership of a regularly airing program from one week to the next, if advertisements are also aired in between those weekly episodes.  So if the first and third airing of the advertisement took place one week apart in successive episodes, and a second airing took place in between, one would expect a large number of repeat viewers between the first and third airing, but not between the second airing and either the first or third airing.

\subsection{Incorporation of Reach into Schedule Optimization}

The approximate reach estimate developed in Subsection~\ref{sec:reach_predict} can be expressed as a polynomial function of the schedule vector $ \Xcia$:
\begin{equation*}
R(\Xcia) = \sum_{i=1}^{N_c} \Xcia S_i \prod_{i'<i} (1-P_{i',i}X_{c,i'}^{(a)}),
\end{equation*}
where $N_c $ is the number of slots available on the channel $c$ under consideration.
Constraints involving reach (or frequency) would become smooth nonlinear constraints, and after relaxation from the integer constraint, could be
approached by the alternating direction method of multipliers \cite{Bazaraa2013}.  

\subsection{Uncertainty Estimation for Reach and Frequency}
For the purpose of building in safety margins in a schedule to avoid disappointing an important advertising client, we might be interested in characterizing the risk that a particular advertising campaign might miss the targets set by an advertiser's bid.   The simplest characterization of uncertainty would be a standard deviation.  If the number of impressions $ I_i $ and the 2-slot characterization of repeat viewership, $ P_{i,i^\prime} $ are directly estimated from historical data by one of the methods described in Section
\ref{sec:predict}, and then those methods could also be used to produce uncertainty estimates.  (Kalman filtering does this automatically.)  The reach and frequency are somewhat complicated functions of these variables, so  in what follows, we describe one simple way we might translate the uncertainty estimates of these variables to an uncertainty estimate for reach and frequency.

We begin by assuming the uncertainty in the estimates of $ \{S_i\}_{i=1}^{N_c} $ and $ \{P_{i,i^\prime}\}_{1 \leq i < i^\prime \leq N_c} $ are all independent, and indicate the mean of a random variable $ Y $ as $ \bar{Y} $ and its standard deviation as $ \sigma(Y) $ (so variance is $ \sigma^2_Y $).  Because the variance of a sum of independent random variables is the sum of the variances of each term, we can therefore express the variance of the reach as a sum of the variances of the new impressions:
\begin{equation*}
\sigma^2(R) = \sum_{j=1}^{N^{(a)}} \sigma^2(S^{\sharp}_{i_j}).
\end{equation*}
The independence assumption does allow the variance of the new impressions, $\sigma^2(S^{\sharp}_{i_a})$ to be worked out in a closed form expression in terms of the mean and standard deviations of $ \{S_i\}_{i=1}^{N_c} $ and $ \{P_{i,i^\prime}\}_{1 \leq i < i^\prime \leq N_c} $, but this expression is quite long.  We therefore compute an approximation to the variance that is valid when the standard deviations of all the constituent random variables are small compared to their means:
\begin{align*}
\sigma (S_i) & \ll \bar{S}_i, \qquad 1 \leq i \leq N_c, \\
\sigma (P_{i,i^\prime})& \ll \bar{P}_{i,i^\prime}, \qquad 1 \leq i < i^\prime \leq N_c.
\end{align*}
Then one can conduct a small noise expansion by writing every random variable in the form $ Y = \bar{Y} + \tilde{Y} $, taking a Taylor expansion up to first order in the fluctuations $ \tilde{S}_i $ and $ \tilde{P}_{i,i^\prime} $, and then computing the variance.  We can thereby obtain:
\begin{equation*}
\sigma^2(R) \approx \sum_{j=1}^{N^{(a)}} \sigma^2 (S_{i_a}) \prod_{j^\prime<j} (1 -\bar{P}_{i_{j^\prime},i_j})^2
+ \sum_{j=1}^{N^{(a)}} \sum_{j^{\prime\prime}=1}^{j-1} \bar{S}_{i_j}^2 \sigma^2 (P_{i_{j^\prime},i_j})\prod_{j^\prime<j,j^\prime \neq j^{\prime\prime}} (1 -\bar{P}_{i_{j^\prime},i_j})^2
\end{equation*}
Actually this small noise expansion can be readily generalized to allow correlations between the random variable models for 
$ \{S_i\}_{i=1}^N $ and $ \{P_{i,i^\prime}\}_{1 \leq i < i^\prime \leq N} $; the same strategy would produce further sums involving the covariances between all pairs of these variables.

Applying the same small noise approximation to the frequency, we obtain an estimate for its standard deviation:
\begin{equation*}
\sigma^2(F) \approx \frac{\sigma^2(S)}{\bar{R}^2} + \frac{\sigma^2(R) \bar{S}^2}{\bar{R}^4}
\end{equation*}
where
\begin{equation*}
\sigma^2(S) =\sum_{j=1}^{N^{(a)}} \sigma^2(S_{i_j}).
\end{equation*}
Similar estimates can be made for the predictions for future numbers of impressions thereby making it possible to estimate the inherent risk in any given schedule.

\section{Conclusions}
\label{sec:conc}
In this report, we analyzed the problem of optimally scheduling advertisements using several different methods. First, we analyzed historical data to obtain trends in the viewership. We found that the viewership was strongly periodic and that deviations from the periodic signal (noise) were approximately bell-shaped.  We supplemented these analyses with predictions from several machine learning algorithms for viewership, from a Bayesian procedure for predicting new program impressions from the program's target demographic, and a measure for comparing programs in order to fill in missing or unknown data. Second, we developed an algorithm, based on binary integer programming, to schedule advertisements. Given orders in the form of a budget, number of impressions desired and demographic targets, the algorithm produces a binary matrix that tells the media company how to schedule advertisements in such a way as to maximize revenue. The algorithm
can be initialized with a schedule generated by a greedy heuristic. Finally, we developed a theoretical framework to quickly estimate the reach (number of new impressions made) of an advertisement. This framework approximates the number of new viewers through historical impressions data and a two-slot function, which gives the fraction of viewers who watched the same advertisement in two time slots.

In summary, mathematical analysis can be an extremely useful tool for understanding how to best schedule advertisements. Techniques from probability, statistics, data science, signal analysis and linear/non-linear programming can all be used to improve and optimize advertising campaigns, give insight into
viewership trends and predict the reach of future television programs.

\begin{backmatter}

\section*{Competing interests}
The authors declare that they have no competing interests.

\section*{Author's contributions}
GS Bhatt, S Burhoe, M Capps, CJ Edholm, S-L Estock, P-W Fok, N Gold, M Houser, P Kramer, H-W Lee, L Rossi, D Shutt, \& VC Yang contributed primarily to the development of methods for predicting viewership.
F El Moustaid,  T Emerson, R Halabi, Q Li, W Li, D Lu, Y Qian, MJ Panaggio,  \& Y Zhou contributed primarily to the formulation and solution of the scheduling optimization problem. 
MJ Panaggio and P-W Fok compiled the content of the manuscript and all contributed to the revision of the manuscript.
\section*{Acknowledgements}
This problem and the data used in this project were provided by Marco Montes de Oca and clypd, Inc.  This work was partially supported by NSF Grant DMS-1261594 through the 2015 Mathematical Problems in Industry Workshop.


\newcommand{\BMCxmlcomment}[1]{}

\BMCxmlcomment{

<refgrp>

<bibl id="B1">
  <title><p>A brief history of advertising in {America}</p></title>
  <aug>
    <au><snm>O'Barr</snm><fnm>WM</fnm></au>
  </aug>
  <source>Advertising \& Society Review</source>
  <publisher>Advertising Educational Foundation</publisher>
  <pubdate>2010</pubdate>
  <volume>11</volume>
  <issue>1</issue>
</bibl>

<bibl id="B2">
  <title><p>Global Media Report 2014</p></title>
  <aug><au><cnm>McKinsey \& Company</cnm></au></aug>
  <pubdate>2014</pubdate>
  <url>http://www.mckinsey.com/~/media/McKinsey/dotcom/client_service/Media
</bibl>

<bibl id="B3">
  <title><p>The Nielsen Company - Solutions: Television</p></title>
  <pubdate>2015</pubdate>
  <url>http://www.nielsen.com/us/en/solutions/measurement/television.html</url>
</bibl>

<bibl id="B4">
  <title><p>Introduction to the theory of statistics.</p></title>
  <aug>
    <au><snm>Mood</snm><fnm>AM</fnm></au>
    <au><snm>Graybill</snm><fnm>FA</fnm></au>
    <au><snm>Boes</snm><fnm>DC</fnm></au>
  </aug>
  <publisher>New York, NY: McGraw-Hill</publisher>
  <pubdate>1974</pubdate>
</bibl>

<bibl id="B5">
  <title><p>A Fresh Look at the {Kalman} Filter</p></title>
  <aug>
    <au><snm>Humpherys</snm><fnm>J</fnm></au>
    <au><snm>Redd</snm><fnm>P</fnm></au>
    <au><snm>West</snm><fnm>J</fnm></au>
  </aug>
  <source>SIAM Review</source>
  <pubdate>2012</pubdate>
  <volume>54</volume>
  <issue>4</issue>
  <fpage>801</fpage>
  <lpage>823</lpage>
</bibl>

<bibl id="B6">
  <title><p>Integer programming</p></title>
  <aug>
    <au><snm>Wolsey</snm><fnm>LA</fnm></au>
  </aug>
  <publisher>New York, NY: John Wiley \& Sons</publisher>
  <pubdate>1998</pubdate>
  <volume>42</volume>
</bibl>

<bibl id="B7">
  <title><p>Integer and combinatorial optimization</p></title>
  <aug>
    <au><snm>Wolsey</snm><fnm>LA</fnm></au>
    <au><snm>Nemhauser</snm><fnm>GL</fnm></au>
  </aug>
  <publisher>Hoboken, NJ: John Wiley \& Sons</publisher>
  <pubdate>2014</pubdate>
</bibl>

<bibl id="B8">
  <title><p>Exploring relaxation induced neighborhoods to improve {MIP}
  solutions</p></title>
  <aug>
    <au><snm>Danna</snm><fnm>E</fnm></au>
    <au><snm>Rothberg</snm><fnm>E</fnm></au>
    <au><snm>Le Pape</snm><fnm>C</fnm></au>
  </aug>
  <source>Mathematical Programming</source>
  <publisher>Springer</publisher>
  <pubdate>2005</pubdate>
  <volume>102</volume>
  <issue>1</issue>
  <fpage>71</fpage>
  <lpage>-90</lpage>
</bibl>

<bibl id="B9">
  <title><p>Preprocessing and probing techniques for mixed integer programming
  problems</p></title>
  <aug>
    <au><snm>Savelsbergh</snm><fnm>MW</fnm></au>
  </aug>
  <source>ORSA Journal on Computing</source>
  <publisher>INFORMS</publisher>
  <pubdate>1994</pubdate>
  <volume>6</volume>
  <issue>4</issue>
  <fpage>445</fpage>
  <lpage>-454</lpage>
</bibl>

<bibl id="B10">
  <title><p>Nonlinear programming: theory and algorithms</p></title>
  <aug>
    <au><snm>Bazaraa</snm><fnm>MS</fnm></au>
    <au><snm>Sherali</snm><fnm>HD</fnm></au>
    <au><snm>Shetty</snm><fnm>CM</fnm></au>
  </aug>
  <publisher>Hoboken, NJ: John Wiley \& Sons</publisher>
  <pubdate>2013</pubdate>
</bibl>

</refgrp>
} 

\end{backmatter}
\end{document}